\documentclass{gtart}

\newcommand{\pathtotrunk}{./}

\usepackage{amsmath,amssymb,amsfonts}
\usepackage{ifpdf}

\ifpdf
\usepackage[pdftex,all,color]{xy}
\else
\usepackage[all,color]{xy}
\fi

\SelectTips{cm}{}

\ifpdf
\newcommand{\pathtodiagrams}{\pathtotrunk diagrams/pdf/}
\else
\newcommand{\pathtodiagrams}{\pathtotrunk diagrams/eps/}
\fi

\newcommand{\mathfig}[2]{{\hspace{-3pt}\begin{array}{c}%
  \raisebox{-2.5pt}{\includegraphics[width=#1\textwidth]{\pathtodiagrams #2}}%
\end{array}\hspace{-3pt}}}

\newcommand{\rotatemathfig}[3]{{\hspace{-3pt}\begin{array}{c}%
  \raisebox{-2.5pt}{\rotatebox{#2}{\includegraphics[height=#1\textwidth]{\pathtodiagrams #3}}}%
\end{array}\hspace{-3pt}}}

\ifpdf
\usepackage[pdftex]{hyperref}
\else
\usepackage[dvips]{hyperref}
\fi
\newcommand{\arxiv}[1]{\href{http://arxiv.org/abs/#1}{\nolinkurl{#1}}}

\theoremstyle{plain}

\newtheorem{prop}{Proposition}[section]
\newtheorem{conj}[prop]{Conjecture}
\newtheorem{thm}[prop]{Theorem}
\newtheorem{lem}[prop]{Lemma}
\newtheorem{cor}[prop]{Corollary}
\newtheorem*{cor*}{Corollary}

\newtheorem{defn}[prop]{Definition}         
\newtheorem*{defn*}{Definition}             

\newenvironment{rem}{\noindent\textsl{Remark.}}{}  
\numberwithin{equation}{section}

\newcounter{comment}
\newcommand{\noop}[1]{}

\def\clap#1{\hbox to 0pt{\hss#1\hss}}
\def\mathllap{\mathpalette\mathllapinternal}
\def\mathrlap{\mathpalette\mathrlapinternal}

\def\mathllapinternal#1#2{%
\llap{$\mathsurround=0pt#1{#2}$}}
\def\mathrlapinternal#1#2{%
\rlap{$\mathsurround=0pt#1{#2}$}}

\newcommand{\Integer}{\mathbb Z}
\newcommand{\Rational}{\mathbb Q}
\newcommand{\Real}{\mathbb R}

\newcommand{\Id}{\boldsymbol{1}}
\renewcommand{\imath}{\mathfrak{i}}
\renewcommand{\jmath}{\mathfrak{j}}

\newcommand{\qRing}{\Integer[q,q^{-1}]}

\newcommand{\To}{\rightarrow}

\newcommand{\Iso}{\cong}

\newcommand{\htpy}{\simeq}

\newcommand{\cone}[3]{C\left(#1 \xrightarrow{#2} #3\right)}
\newcommand{\pairing}[2]{\left\langle#1 ,#2 \right\rangle}

\newcommand{\card}[1]{\sharp{#1}}

\newcommand{\bdy}{\partial}
\newcommand{\compose}{\circ}

\newcommand{\Cat}{\mathcal{C}}

\newcommand{\psmallmatrix}[1]{\left(\begin{smallmatrix} #1 \end{smallmatrix}\right)}

\newcommand{\qdim}{\operatorname{dim_q}}

\newcommand{\directSum}{\oplus}
\newcommand{\DirectSum}{\bigoplus}
\newcommand{\tensor}{\otimes}
\newcommand{\Tensor}{\bigotimes}

\newcommand{\su}[1]{\mathfrak{su}_{#1}}
\newcommand{\csl}[1]{\mathfrak{sl}_{#1}}
\newcommand{\uqsl}[1]{U_q\left(\csl{#1}\right)}

\newcommand{\Cobl}{{\mathcal Cob}_{/l}}
\newcommand{\Cob}[1]{{\mathcal Cob}\left(\su{#1}\right)}

\newcommand{\Mat}[1]{\mathbf{Mat}\left(#1\right)}
\newcommand{\Kar}[1]{\mathbf{Kar}\left(#1\right)}
\newcommand{\Inv}[1]{\operatorname{Inv}\left(#1\right)}
\newcommand{\Hom}[3]{\operatorname{Hom}_{#1}\left(#2,#3\right)}
\newcommand{\End}[1]{\operatorname{End}\left(#1\right)}

\newcommand{\Gr}[2]{\text{Gr}(#1 \subset #2)}
\newcommand{\Flag}[3]{\text{Flag}(#1 \subset #2 \subset #3)}

\def\llbracket{\left[\!\!\left[}
\def\rrbracket{\right]\!\!\right]}
\newcommand{\Kh}[1]{\llbracket#1\rrbracket}

\newcommand{\TL}{\mathcal{TL}}
\newcommand{\Foam}[1]{\mathbf{Foam}\left(#1\right)}

\newcommand{\directSumStack}[2]{{\begin{matrix}#1 \\ \DirectSum \\#2\end{matrix}}}
\newcommand{\directSumStackThree}[3]{{\begin{matrix}#1 \\ \DirectSum \\#2 \\ \DirectSum \\#3\end{matrix}}}

\newcommand{\grading}[1]{{\color{blue}\{#1\}}}
\newcommand{\shift}[1]{\left[#1\right]}
\newcommand{\pa}{\mathcal{P}}

\usepackage{breakurl}

\ifpdf
\usepackage[pdftex]{graphicx}
\else
\usepackage[dvips]{graphicx}
\fi

\usepackage{color}

\title{On Khovanov's cobordism theory for $\su{3}$ knot homology.}

\author{Kim~Morrison}
\address{
   Department of Mathematics\\
   University of California, Berkeley\\
   Berkeley CA 94720\\
   USA
} \email{scott@math.berkeley.edu}
\urladdr{http://math.berkeley.edu/\char126 scott}

\author{Ari~Nieh}
\email{ari@math.berkeley.edu}
\urladdr{http://math.berkeley.edu/\char126 ari}

\date{
  First edition: December 24, 2006.
  This edition: \today.
}


\begin{document}

\begin{abstract}
We reconsider the $\su{3}$ link homology theory defined by Khovanov
in \cite{MR2100691}  
and generalized by Mackaay and Vaz in \cite{math.GT/0603307}. With
some slight modifications, we describe the theory as a map from the
planar algebra of tangles to a planar algebra of (complexes of) `cobordisms with
seams' (actually, a `canopolis'), making it local in the sense of
Bar-Natan's local $\su{2}$ theory of \cite{MR2174270}. 

We show that this `seamed cobordism canopolis' decategorifies to
give precisely what you'd both hope for and expect: Kuperberg's
$\su{3}$ spider defined in \cite{MR1403861}. 
We conjecture an answer to an even more interesting question about the decategorification of the Karoubi envelope of our cobordism theory.

Finally, we describe how the theory is actually completely
computable, and give a detailed calculation of the $\su{3}$ homology
of the $(2,n)$ torus knots.
\end{abstract}

\maketitle

\tableofcontents

\section{Introduction}
Bar-Natan formulated a highly geometric version of Khovanov homology
in \cite{MR2174270}.  His approach uses the language of planar
algebras for the construction of the complex.  In particular, it has
the pleasant feature of being a local theory, which makes it useful
for fast `divide and conquer' computations \cite{math.GT/0606318}.

Khovanov constructed a homology theory of links that
categorifies the $\su{3}$ quantum knot invariant in \cite{MR2100691}.  Mackaay and Vaz generalized this theory in \cite{math.GT/0603307}.  In the spirit of
Bar-Natan, we provide a local perspective on this knot homology. Our
formulation uses a planar algebra of categories (a `canopolis') as
the setting for the complex.

The $\su{3}$ quantum knot invariant is determined by the following
formulas, which should be thought of as a map of planar algebras:
\begin{align*}
 \mathfig{0.06}{webs/positive_crossing}&  \mapsto  \phantom{- q^{-3}} \mathllap{q^2}  \mathfig{0.06}{webs/two_strand_identity} - \phantom{q^{-2}} \mathllap{ q^3} \mathfig{0.06}{webs/upwards_H_diagram} \\
 \mathfig{0.06}{webs/negative_crossing}&  \mapsto - q^{-3} \mathfig{0.06}{webs/upwards_H_diagram} + q^{-2} \mathfig{0.06}{webs/two_strand_identity}
\end{align*}

This sends an oriented link diagram to a $\qRing$-linear
combination of oriented planar graphs with trivalent vertices (`webs').  We then
evaluate these webs using the relations of Kuperberg's $\su{3}$
spider \cite{MR1403861}
\begin{align}\label{spider-relations}
   \mathfig{0.045}{webs/clockwise_circle} & = q^2 + 1 + q^{-2}\\
   \mathfig{0.045}{webs/bubble} & = q\mathfig{0.009}{webs/tall_strand} + q^{-1}\mathfig{0.009}{webs/tall_strand} \\
   \mathfig{0.08}{webs/oriented_square}
     & = \mathfig{0.08}{webs/two_strands_horizontal} + \mathfig{0.08}{webs/two_strands_vertical}
\end{align}
to obtain a polynomial invariant of links.

Just as the categorified version of (one variation of) the Kauffman skein relation for
the Jones polynomial\footnote{This isn't quite the quantum $\su{2}$
skein theory; see \cite{morrison-walker}.}
\begin{align*}
 \mathfig{0.08}{webs/positive_crossing}&  \mapsto  q \mathfig{0.08}{webs/two_unoriented_strands_vertical} - q^2 \mathfig{0.08}{webs/two_unoriented_strands_horizontal}
\end{align*}
becomes the following complex in Khovanov's theory,
\begin{equation*}
\xymatrix@R-1mm{
 \mathfig{0.06}{webs/positive_crossing} \ar@{|->}[r] & \Bigg( \bullet \ar[r]
    & q \mathfig{0.06}{webs/two_unoriented_strands_vertical} \ar[r]^{\mathfig{0.05}{cobordisms/unoriented_saddle1}} & q^2 \mathfig{0.06}{webs/two_unoriented_strands_horizontal} \ar[r] & \bullet \Bigg) \\
}
\end{equation*}
we should expect the categorified $\su{3}$
invariant to associate to a crossing some two step complex, with something like a cobordism for the differential.
However, since the diagrams in the $\su{3}$ spider have singularities, the
category of cobordisms can't suffice; therefore, we'll work with
seamed cobordisms (or `foams') that allow singular seams where three
half-planes meet:
\begin{equation*}%
\xymatrix@R-1mm{
 \mathfig{0.06}{webs/positive_crossing} \ar@{|->}[r] & & \Bigg( \bullet \ar[r]
    & q^2 \mathfig{0.06}{webs/two_strand_identity} \ar[r]^{\mathfig{0.06}{cobordisms/zip}} & q^3 \mathfig{0.06}{webs/upwards_H_diagram} \ar[r] & \bullet \Bigg) \\
 \mathfig{0.06}{webs/negative_crossing} \ar@{|->}[r] & \Bigg( \bullet \ar[r] &
    q^{-3} \mathfig{0.06}{webs/upwards_H_diagram} \ar[r]^{\mathfig{0.06}{cobordisms/unzip}} & q^{-2} \mathfig{0.06}{webs/two_strand_identity} \ar[r] & \bullet \Bigg) &
}%
\end{equation*}

We'll describe this construction in detail, essentially paralleling
the work of Khovanov and of Mackaay and Vaz, with some minor
differences which we find appealing.\footnote{Much of our work was done before the appearance of
\cite{math.GT/0603307}, which perhaps partially excuses our giving a
self-contained development of the theory.} For most of the paper, it
isn't necessary to have read their work (although \S
\ref{sec:mackaay-vaz} which explicitly compares the details of our
construction with that of Khovanov and of Mackaay and Vaz assumes this). We
emphasize the local nature of our construction, giving automatic
proofs of Reidemeister invariance, following Bar-Natan's
simplification algorithm, in \S \ref{ssec:isotopy-invariance}.
Later, in \S \ref{sec:2-n-torus-knots}, we provide explicit detailed
calculations of the $\su{3}$ Khovanov invariant for the $(2,n)$
torus knots.

Our version of this invariant associates to every tangle an
up-to-homotopy complex in the canopolis of foams.  In \S
\ref{sec:decategorification}, we prove `decategorification'
results both for this canopolis and for Bar-Natan's canopolis of
cobordisms corresponding to the original Khovanov homology.  Roughly
speaking, this involves collapsing the categorical structure of the
canopolis (taking the split Grothendieck group) while preserving
its planar algebra structure.  The decategorification of Bar-Natan's
canopolis is the Temperley-Lieb planar algebra. Similarly, the
decategorification of the canopolis of foams is the Kuperberg's
$\su{3}$ spider.  As we will see, the $\su{3}$ case requires more
complicated techniques, because the morphisms are much harder to
classify than the cobordisms in the $\su{2}$ canopolis.  Among these
techniques is a kind of duality: in \S
\ref{ssec:su_3-decategorification} we'll produce isomorphisms
$\Hom{}{U \otimes V}{W} \cong \Hom{}{U}{W \otimes V^*}$ in the
canopolis of $\su{3}$ foams, which we think of as meaning
that it's secretly a `spatial algebra' (i.e. a higher dimensional analogue of a planar algebra), not just a canopolis.

Some interesting things happen in the $\su{3}$ theory which have no
analogues for $\su{2}$. In particular, there are grading $0$
morphisms other than the identity between irreducible diagrams.
We'll discuss an example in which the identity morphism can be
written as a sum of orthogonal idempotents, and make a conjecture
about the decategorification of the Karoubi envelope. (The Karoubi
envelope is the category we get by adding in all idempotents as
extra objects.) A further conjecture says that the minimal
idempotents correspond to the dual canonical basis in the $\su{3}$
spider \cite{MR1680395}.

\section{Preliminaries}

\subsection{Locality, or, ``What is a planar algebra?''}
A planar algebra is a gadget specifying how to combine objects in
planar ways. They were introduced in \cite{math.QA/9909027} to study
subfactors, and have since found more general use.

In the simplest version, a planar algebra $\pa$ associates a vector
space $\pa_k$ to each natural number $k$ (thought of as a disc in
the plane with $k$ points on its boundary) and a linear map
$\pa(T)  : \pa_{k_1} \tensor \pa_{k_2} \tensor \cdots \tensor
\pa_{k_r} \To \pa_{k_0}$ to each `spaghetti and meatballs' diagram
$T$, for example
$$\mathfig{0.2}{planar_tangles/cubic_tangle_example},$$ with internal
discs with $k_1, k_2, \ldots, k_r$ points, and $k_0$
points on the external disc. These maps (the 'planar operations')
must satisfy certain properties: radial spaghetti induce
identity maps, and composition of the maps $\pa(T)$ is compatible
with the obvious composition of spaghetti and meatballs diagrams by
gluing one inside the other.

For the exact details, which are somewhat technical, see
\cite{math.QA/9909027}.

{
\newcommand{\labels}{\mathfrak{L}}
Planar algebras also come in more subtle flavors. Firstly, we can
introduce a label set, and associate a vector space to
each disc with boundary points colored by this label set. (The
simplest version discussed above thus has a singleton label set, and
the discs are indexed by the number of boundary points.) The planar
tangles must now have arcs colored using the color set, and the
rules for composition of diagrams require that labels match up.
We can also have a oriented label set; the label set has an involution and the arcs carry both an orientation and a label, modulo reversing both.
Secondly, we needn't restrict ourselves to vector spaces and linear
maps between them; a planar algebra can be defined over an arbitrary
monoidal category, associating objects to discs, and morphisms to
planar tangles. Thus we might say ``$\mathcal{P}$ is a planar
algebra over the category $\Cat$ with label set $\labels$.''
\footnote{A subfactor planar algebra is defined over
$\operatorname{Vect}$, and has a 2 element label set. One imposes an
additional condition that only discs with an even number of boundary
points and with alternating labels have non-trivial vector spaces
attached. There is also a positivity
condition. See \cite[\S4]{math.GT/0606318}.}%
}

A `canopolis', introduced by Bar-Natan in
\cite{MR2174270}\footnote{He called it a `canopoly', instead, but
we're taking the liberty of fixing the name here.}\footnote{See also
\cite{math.GT/0610650} for a description of Khovanov-Rozansky
homology \cite{math.QA/0401268, math.QA/0505056} using
canopolises.}, is simply a planar algebra defined over some category
of categories, with monoidal structure given by cartesian product.
Thus to each disc, we associate some category of a specified type. A
planar tangle then induces a functor from the product of internal
disc categories to the outer disc category, thus taking a tuple of
internal disc objects to an external disc object, and a tuple of
internal disc morphisms to an external disc morphism. It is
picturesque to think of the objects living on discs, and the
morphisms in cans, whose bottom and top surfaces correspond to the
source and target objects. Composition of morphisms is achieved by
stacking cans vertically, and the planar operations put cans side by
side.

The functoriality of the planar algebra operations ensure that we
can build a `city of cans' (hence the name canopolis) any way we
like, obtaining the same result: either constructing several towers
of cans by composing morphisms, then combining them horizontally,
or constructing each layer by combining the levels of all the towers
using the planar operations, and then stacking the levels
vertically.

\subsection{The $\su{2}$ cobordism theory}
We will now briefly recall the canopolis defined by Bar-Natan in
\cite{MR2174270}, and used in his local link homology theory.

Slightly modifying Bar-Natan's notation, $\Cob{2}$ is our name for
his $\Cobl^3$, the canopolis of cobordisms in cans modulo the
$\su{2}$ relations.

The objects of $\Cob{2}$ consist of planar tangle diagrams:
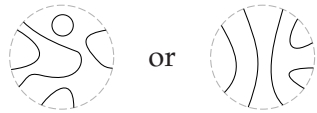
\begin{equation*}
\mathfig{0.1}{planar_tangles/TL_example1} \quad \text{or} \quad
\mathfig{0.1}{planar_tangles/TL_example2}
\end{equation*}

equipped with the obvious planar algebra structure\footnote{We may
think of this as the free planar algebra with no generators.}.

Let $R_0$ be any commutative ring in which 2 is invertible.  If
$D_1$ and $D_2$ are diagrams with identical boundary, a morphism
between them is a formal $R_0$-linear combination of cobordisms from
$D_1$ to $D_2$ modulo the following local relations:
\begin{align*}%
\mathfig{0.1}{cobordisms/sphere} & = 0 &
\mathfig{0.1}{cobordisms/torus} & = 2 \end{align*}
\begin{align}\label{su_2-neck-cutting}
\rotatemathfig{0.05}{90}{cobordisms/cylinder} & = \frac{1}{2}
\rotatemathfig{0.05}{90}{cobordisms/neck_cutting_left} + \frac{1}{2}
      \rotatemathfig{0.05}{90}{cobordisms/neck_cutting_right}
\end{align}
The planar algebra structure on morphisms is given by plugging
cans into $T \times [0,1]$, where $T$ is a spaghetti and
meatballs diagram, as in this example:

$$\mathfig{0.2}{canopolis/tangle_action_on_cans}.$$

We refine the theory by introducing a grading on the canopolis.  We
equip the objects of $\Cob{2}$ with a formal grading shift, so that
they are of the form $q^m D$, where $m$ is an integer\footnote{This
is Bar-Natan's $D\grading{m}$.}. (We will, however, sometimes
suppress the grading for simplicity, or conflate diagrams with
objects when it is convenient.)  We let grading shifts add under
planar algebra operations. The degree of a cobordism $C$ from
$q^{m_1}D_1$ to $q^{m_2}D_2$ is defined as $\chi(C) - k/2 + m_2 -
m_1$, where $\chi$ is the Euler characteristic and $k$ is the number
of boundary points of $D_i$. It is not hard to see that degrees are
additive under both composition and planar
operations.\footnote{Observe that $\chi(c)-\frac{k}{2}$ and $m_2 -
m_1$ are additive separately.} Note also that the local relations
are degree-homogeneous, and therefore this grading descends to the
quotient.

We can further introduce formal direct sums, and allow matrices of
morphisms between direct summands. This is the matrix category
construction, applied to each category in our canopolis. We denote
the result $\Mat{\Cob2}$.

\subsubsection{The structure of morphisms in $\Cob{2}$}
The structure of this canopolis has been thoroughly analyzed
elsewhere, in Dror's paper \cite[\S 9]{MR2174270} and in Gad Naot's
\cite{math.GT/0603347}. We will need one of their results.

First, note that almost all closed surfaces in $\Cob{2}$ can be
evaluated as scalars.  In fact, applying the `neck-cutting'
relation (\ref{su_2-neck-cutting}) shows that they are all zero except for the surfaces of
genus one and three.  We saw above that the torus was equal to 2,
but there is no {\it a priori} way to evaluate the surface of genus
three.  Therefore, we absorb it into our ground ring, letting $R =
R_0[\mathfig{0.1}{cobordisms/triple_torus}]$.

\begin{prop}\label{thm:su_2-classification}
For any two diagrams $D_1$ and $D_2$, let $l$ be the number of
components of $D_1 \cup D_2 \cup (\partial \times [0,1])$.  Consider
the set of cobordisms $C \in \Hom{\Cob{2}}{D_1}{D_2}$ such that
every component of $C$ is either a disc or a punctured torus.  These
cobordisms form a basis for $\Hom{\Cob{2}}{D_1}{D_2}$ over $R$.
\end{prop}

Note that such cobordisms must have exactly $l$ components, and the
boundary of each component is a single component of $D_1 \cup D_2
\cup (\partial \times [0,1])$.

\begin{rem}
This classification requires the neck-cutting relation, and only
holds when $2$ is invertible.  (See \cite{math.GT/0603347} for
details otherwise.)
\end{rem}

We call a diagram `non-elliptic' \footnote{This is the obvious
extension of Kuperberg's meaning of `non-elliptic' in \cite{MR1403861} to the
$\su{2}$ case.} if it contains no circles. By the previous result
and some Euler characteristic calculations, we get:

\begin{cor}\label{cor:endomorphisms-grading}
Endomorphisms of a non-elliptic diagram are all in non-positive
degree.
\end{cor}

\begin{cor}\label{cor:only-identity}
If a nonzero endomorphism of a non-elliptic diagram factors through
a different non-elliptic diagram, then it necessarily has negative
grading.
\end{cor}

\begin{rem}
It's easy to see that elliptic diagrams have positively graded
endomorphisms; for example, a circle which dies and is born again,
each time via a disc cobordism, has grading $+2$.
\end{rem}

This classification also yields a description of the `sheet
algebra' for the $\su{2}$ canopolis:

\begin{cor}\label{cor:su_2-sheet-algebra}
Let S be the diagram consisting of a single arc. Then
$$\End{S} =
R\left[\left.\mathfig{0.03}{cobordisms/sheet_algebra/handle}\right]
\right/ \left \langle
\mathfig{0.03}{cobordisms/sheet_algebra/handle}^2 - \frac{1}{2}
\mathfig{0.025}{cobordisms/sheet_algebra/sheet} \! \! \!
\mathfig{0.1}{cobordisms/triple_torus} \right \rangle.$$
\end{cor}

\section{The $\su{3}$ cobordism theory}

\subsection{Seamed cobordisms, and the $\su{3}$ theory} We now describe
$\Cob3$, the analogous canopolis of `seamed cobordisms' associated to
$\su{3}$. The objects consist of `webs' -- elements of the planar
algebra freely generated by the trivalent vertices
$$\mathfig{0.10}{webs/inwards_vertex} \qquad \text{and}
\qquad \mathfig{0.10}{webs/outwards_vertex}.$$%
(It's a planar algebra whose label set consists of just two labels:
`in' and `out'.)
Let $S$ be a commutative ring in which 2 and 3
are invertible.  The set of morphisms between two webs with the same
boundary will be an $S$-module generated by `seamed cobordisms',
also called `foams'.

The local model for a seamed cobordism is the space $Y \times [0,1]$,
the space obtained by gluing together three copies of $[0,1] \times [0,1]$
along $[0,1] \times \{0\}$, with orientations on the three squares, all
inducing the same orientation on the common $[0,1] \times \{0\}$, along
with a cyclic orientation of the three squares.\footnote{We say that
a seamed cobordism $C$ is locally modeled on $Y \times [0,1]$ in the
same sense that that a topological $n$-manifold is modeled on
(topological) $\Real^n$. We mean that for every point $p$ of $C$,
there is a point $p'$ of $Y$, neighborhoods $p \in U_p \subset C$
and $p' \in U'_{p'} \subset Y \times [0,1]$ and a  bijection $f_p : U_p
\To U'_{p'}$. Moreover, the `transition maps' $f_p^{-1} f_q$ should
preserve the local structure specified for $Y \times [0,1]$; in
particular, the topological structure and, more importantly, the
orientation data.}

\begin{defn}
Given two webs, $D_1$ and $D_2$, drawn in a disc, both with boundary
$\partial$, a seamed cobordism from $D_1$ to $D_2$ is a
2-dimensional CW-complex \footnote{We don't care about the actual cell
decomposition, of course.} $F$ (the `foam') with
\begin{itemize}
\item exactly three $2$-cells meeting along each singular $1$-cell,
\item a cyclic ordering on those three $2$-cells,
\item orientations on the $2$-cells, compatible with the cyclic
orderings,
\item and an identification of the boundary of $F$ with $D_1 \cup D_2
\cup (\partial \times [0,1])$ such that
\begin{itemize}
\item the orientations on the sheets induce the orientations on the
edges of $D_1$, and the opposite orientations on the edges of $D_2$,
\item and the cyclic orderings around the singular seams agree with the
cyclic orderings around a vertex in $D_1$ or $D_2$ given by its
embedding in the disc; the anticlockwise ordering for `inwards' vertices, the clockwise ordering for `outwards' vertices.
\end{itemize}
\end{itemize}
\end{defn}
We think of such a foam as living inside the `can' $D^2
\times [0,1]$, even though it is not embedded there; there's just an
identification of its boundary with a subset of the surface of the
can.

Compositions, both vertical (everyday composition of morphisms in a
category) and horizontal (the action of planar tangles on
morphisms), are almost trivial to describe. To compose vertically,
we stack cans on top of each other, and to compose horizontally using
a spaghetti and meatballs diagram $T$, we glue together $T \times [0,1]$ with the input
cans.

As before, to put a grading on our canopolis, we endow diagrams with
formal grading shifts written as factors of $q$.  The degree of a
cobordism $C$ from $q^{m_1} D_1$ to $q^{m_2} D_2$ is defined as
\begin{equation}
\label{eq:su_3-grading}
 \deg{C} = 2\chi(C) - \card{\bdy} + \frac{\card{V}}{2} + m_2 - m_1,
\end{equation}
where $\card{\bdy}$ is the number of boundary points of $D_i$ and $\card{V}$ is the
total number of trivalent vertices in $D_1$ and $D_2$.  We leave it
to the reader to check that this is additive under canopolis
operations.

It is not hard to verify that this canopolis of $\su{3}$ foams is
generated (as a canopolis!) by the morphisms cup, cap, saddle, zip,
and unzip (after \cite{MR2048108}):
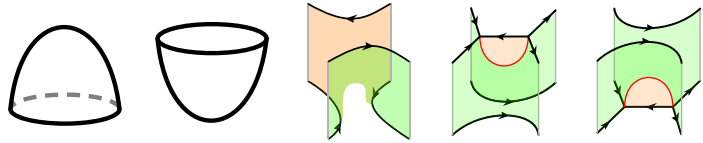
\begin{equation*}
\rotatebox[origin=c]{90}{$\mathfig{0.1}{cobordisms/cap_bdy_left}$}
\quad
\rotatebox[origin=c]{90}{$\mathfig{0.1}{cobordisms/cap_bdy_right}$}
\quad \mathfig{0.1}{cobordisms/saddle1} \quad
\mathfig{0.1}{cobordisms/zip} \quad \mathfig{0.1}{cobordisms/unzip}
\end{equation*}

As a little piece of nomenclature, we'll introduce the cobordism we
call a `choking torus',
$\mathfig{0.1}{cobordisms/handle_disc_bdy_left}$. Whenever you see
this, you should assume the cyclic ordering at the seam is
`bulk/handle/disc'.

\subsection{Local relations}
\label{sec:local-relations}%
We now introduce local relations on the modules of seamed
cobordisms. These are motivated in two ways:
\begin{enumerate}
\item We expect that the canopolis of seamed cobordisms should have
isomorphisms reflecting the relations appearing in the $\su{3}$
spider.
\item We intend to construct an invariant of tangles, valued
in complexes of seamed cobordisms.
\end{enumerate}

We'll see both of these motivations validated, in sections
\S\ref{ssec:isomorphisms} and \S\ref{ssec:simplification}
respectively.

\begin{itemize}
\item`Closed foam' relations:
\begin{align}
\label{eq:closed_foams}%
\mathfig{0.1}{cobordisms/sphere} & = 0 &
\mathfig{0.1}{cobordisms/torus} & = 3 \\
\notag
\mathfig{0.2}{cobordisms/double_torus} & = 0 &
\mathfig{0.24}{cobordisms/triple_torus_disc} & = 0
\end{align}

\item The `neck cutting' relation:
\begin{equation}
\label{eq:neck_cutting}
 \rotatemathfig{0.05}{90}{cobordisms/cylinder} = \frac{1}{3}   \rotatemathfig{0.05}{90}{cobordisms/neck_cutting_left}
      - \frac{1}{9} \rotatemathfig{0.05}{90}{cobordisms/neck_cutting_middle}
      + \frac{1}{3} \rotatemathfig{0.05}{90}{cobordisms/neck_cutting_right}
\end{equation}

\item The `airlock' relation:
\begin{equation}
\label{eq:airlock}
    \rotatemathfig{0.055}{90}{cobordisms/airlock} = -
    \rotatemathfig{0.055}{90}{cobordisms/cup_cap}
\end{equation}

\item The `tube' relation
\begin{equation}
\label{eq:tube_relation}
    \mathfig{0.15}{cobordisms/tube_relation/tube} =
    \frac{1}{2}
    \mathfig{0.15}{cobordisms/tube_relation/two_bubbles_lower_kiss} +
    \frac{1}{2}
    \mathfig{0.15}{cobordisms/tube_relation/two_bubbles_upper_kiss}
\end{equation}
The small (green) circles here indicate the two sheets coming together; they're a composition, zip followed by unzip.

\item The `three rocket' relation:
\begin{equation}
\label{eq:rocket_relation}
  \mathfig{0.15}{cobordisms/rocket_relation/rocket_z}
 + \mathfig{0.15}{cobordisms/rocket_relation/rocket_x}
 + \mathfig{0.15}{cobordisms/rocket_relation/rocket_y} = 0
\end{equation}

\item The `seam-swap' relation: reversing the cyclic order of the three 2-cells
attached to a closed singular seam is equivalent to multiplication by -1.
\end{itemize}

As consequences of the above relations, it is not hard to derive the
following:

\begin{itemize}
\item The sheet relations:
\begin{align}
  \mathfig{0.06}{cobordisms/sheet_algebra/blister} & = 0 &
 \label{eq:choking-torus-multiplication}%
   \mathfig{0.075}{cobordisms/sheet_algebra/two_handle_discs} & =
   -3 \mathfig{0.075}{cobordisms/sheet_algebra/handle} \\
 \label{eq:torus-choking-torus}%
 \mathfig{0.075}{cobordisms/sheet_algebra/handle_handle_disc} & = 0 &
  \mathfig{0.075}{cobordisms/sheet_algebra/two_handles} & = 0
\end{align}

The `blister' relation follows directly from
seam-swapping. The `choking torus multiplication' relation on the
first line follows from applying neck-cutting in reverse. The equations in the last line follow from neck cutting, and the closed foam relations.

\item The `bamboo' relation:
\begin{equation}
\label{eq:bamboo}
  \mathfig{0.15}{cobordisms/bamboo/bamboo} =
    \frac{1}{3}
    \mathfig{0.15}{cobordisms/bamboo/bamboo_left} +
    \frac{1}{3}
    \mathfig{0.15}{cobordisms/bamboo/bamboo_right}
\end{equation}

which follows from neck-cutting one side of the bamboo, then
reducing terms via airlocks and blisters.

\end{itemize}

As before, we introduce formal direct sums of the objects and
matrices of morphisms, yielding a canopolis we call $\Mat{\Cob{3}}$.

\subsection{Consistency}
The purpose of this section is two-fold. First, we want to provide a
set of assumptions, plausibly desirable in any categorification of
the $\su{3}$ planar algebra, which allows us to to derive the
relations described in the previous section. Second, we prove the
following result:

\begin{thm}
\label{thm:consistency}%
The local relations of \S \ref{sec:local-relations} are consistent,
in the sense that
$$\Hom{\Cob{3}}{\emptyset}{\emptyset} \neq 0.$$
\end{thm}

These two goals are related. In the process of justifying the local
relations, we will divide them into two classes: the `evaluation
relations', and the `local kernel' relations. The evaluation
relations are the `closed foam' relations, `seam swapping', `neck cutting' and
`airlock'. The `local kernel'
relations are `tube' and `rocket'. We begin by showing the
evaluation relations follow from some appealing assumptions. We then
show that these relations, living up to their name, suffice to
evaluate any closed foam. Further, in \S
\ref{sec:evaluation-consistency} we'll show they're consistent;
denoting the canopolis in which we only impose the evaluation
relations by $\Cob{3}^\text{ev}$, we have
\begin{lem}
\label{lem:evaluation-consistency}%
$$\Hom{\Cob{3}^{\text{ev}}}{\emptyset}{\emptyset} = S.$$
\end{lem}
It's then time to introduce the local kernel relations. The
canopolis $\Cob{3}^{\text{ev}}$ is an unsatisfactory one, in the
sense that it is `degenerate' or has a `local kernel': non-zero
foams with boundary, all of whose completions to a closed foam
are zero. In a slightly different guise, Khovanov proved the following lemma in
\cite{MR2100691}:
\begin{lem}
\label{lem:tube-and-rocket-in-local-kernel}%
The tube relation and rocket relation are in the local kernel
(justifying the name `local kernel relations').
\end{lem}
We'll show in \S \ref{sec:local-kernel} that
\begin{lem}
\label{lem:local-kernel-generators}%
The local kernel is generated,
as a canopolis ideal, by the tube and rocket relations.
\end{lem}

We thus impose the local kernel as additional relations, and
together Lemmas \ref{lem:evaluation-consistency} and
\ref{lem:local-kernel-generators} imply Theorem
\ref{thm:consistency}.

\subsubsection{Explaining the relations}%
\label{sec:explaining-relations}%
We now set out some plausible assumptions one might make about any
categorification of the $\su{3}$ spider. (Perhaps these assumptions
might be useful to someone categorifying something else, as well!)

Firstly, we'll ask, without much motivation, for the grading rule
given previously; the grading of a morphism is given by twice its
Euler characteristic, as in Equation \eqref{eq:su_3-grading}.

We'll just have to pull the `seam-swapping' relation described
earlier out of a hat.\footnote{Note though, that it's the $n=3$
special case of the idea described in \cite[\S 6]{MR2100691} that if
the `$k$-sheets' of an $\su{n}$ foam were to be labeled by elements of
the cohomology ring of $\Gr{k}{n}$, then the relations around a seam
should be the kernel of the map $\Tensor_i H^\bullet(\Gr{k_i}{n})
\To H^\bullet\left(\Flag{k_1}{k_1 + k_2 \subset \dotsb}{\left(\sum_i
k_i\right)}\right)$ induced by the `take orthogonal complements' map
at the geometric level.} This relation kills off certain closed
foams, amongst them the `theta' foam, the `blistered torus'
$\mathfig{0.065}{cobordisms/torus_right_blister}$ (in fact, any foam
with a blister) and $\mathfig{0.11}{cobordisms/double_torus_disc}$.

We'll then put in by hand a few relations motivated by the desire
that $\Hom{\Cob{3}}{\emptyset}{\emptyset}$, the space of closed
foams, as a graded $S$-module, be just $S$ generated by the empty
foam. Later, we'll see that the relations we've imposed do in fact
imply this. First of all, we force the sphere to be zero (it's in
positive degree) and the torus to be some element of $S$. We'll
assume, in fact, that the torus is invertible. Briefly, we'll write
$t$ for this value, but very shortly discover that $t=3$. Further,
various closed foams with negative degrees are forced to be zero,
such as
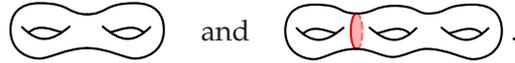
\begin{equation*}
 \mathfig{0.15}{cobordisms/double_torus}  \quad \text{and} \quad
 \mathfig{0.21}{cobordisms/triple_torus_disc}.
\end{equation*}
(However, see \S \ref{sec:relaxing-relations} for a discussion of
the variation in which we just ask that
$\Hom{\Cob{3}}{\emptyset}{\emptyset}_{>0} = 0$ and
$\Hom{\Cob{3}}{\emptyset}{\emptyset}_0$ is 1-dimensional.)

\newcommand{\cc}{\mathfig{0.025}{webs/clockwise_circle}}%
\newcommand{\ac}{\mathfig{0.025}{webs/anticlockwise_circle}}%
Next, we'll ask that $\Hom{\Cob{3}}{\cc}{\emptyset}$ is a free
module of rank $3$, and in fact with graded dimension $q^2 + 1 +
q^{-2}$, on the basis that we expect this graded dimension to agree
with the evaluation of $\cc$ in
the $\su{3}$ spider. Since the cobordisms
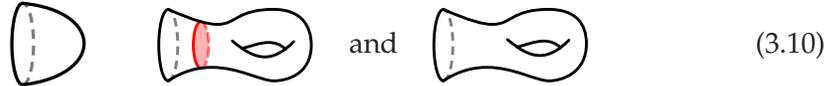
\begin{equation}
\label{eq:disc-morphisms}
 \mathfig{0.075}{cobordisms/cap_bdy_left} \qquad
 \mathfig{0.15}{cobordisms/handle_disc_bdy_left} \quad \text{and}
 \quad
 \mathfig{0.15}{cobordisms/handle_bdy_left}
\end{equation}
lie in this morphism space, with gradings $2$, $0$ and $-2$
respectively, we'll further ask that in fact the morphism space is
freely generated by these three cobordisms. (Unsurprisingly, we'll
ask the same thing for $\ac$.) Remember there are two variations of
the middle cobordism above, differing in the cyclic ordering of the
sheets at the seam; the two cyclic orderings only differ by a sign, however, by the seam-swapping relation.

Further, we'll ask that $\Hom{}{\cc}{\cc} \Iso
\Hom{}{\cc\ac}{\emptyset}$, with the isomorphism given by isotopy.
This behavior will follow from any good notion of duality in a
categorification; moreover, it certainly happens in the $\su{2}$ canopolis,
and we'll see the appropriate generalization to arbitrary diagrams
in \S \ref{sec:local-kernel}. Even more, we'll ask that the obvious map
$\Hom{}{\cc}{\emptyset} \tensor \Hom{}{\ac}{\emptyset} \to
\Hom{}{\cc\ac}{\emptyset}$, given by disjoint union, is actually an
isomorphism; again, we'll later see that this is generally true.

With these relatively benign constraints, we can get a long way!
Firstly, looking at the degree $4$ piece of $\Hom{}{\cc}{\cc}$, we
see it's $1$ dimensional, and so the `airlock'
$\mathfig{0.075}{cobordisms/airlock}$ must be proportional to
$\mathfig{0.075}{cobordisms/cup_cap}$. We'll
declare\footnote{We could try an arbitrary constant here,
$\mathfig{0.055}{cobordisms/airlock} = - \mu
\mathfig{0.055}{cobordisms/cup_cap}$, say. The argument above would
continue much the same, except that we wouldn't be able to find an
analogue of the tube and rocket relations in the local kernel.} that
$$\mathfig{0.09}{cobordisms/airlock} = -
\mathfig{0.09}{cobordisms/cup_cap}.$$

Next, looking at the degree $0$ piece, we see a $3$ dimensional
space. Writing down $4$ obvious cobordisms here,
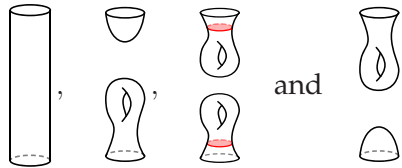
\begin{equation*}
 \rotatemathfig{0.04}{90}{cobordisms/cylinder}, \quad
 \rotatemathfig{0.04}{90}{cobordisms/neck_cutting_left}, \quad
      \rotatemathfig{0.04}{90}{cobordisms/neck_cutting_middle} \quad
      \text{and} \quad \rotatemathfig{0.04}{90}{cobordisms/neck_cutting_right}
\end{equation*}
we see there must be some relation amongst them (this will turn out
to be neck cutting, of course), which we'll suppose is of the form
\begin{equation*}
 \rotatemathfig{0.04}{90}{cobordisms/cylinder} =
 x \rotatemathfig{0.04}{90}{cobordisms/neck_cutting_left}+
 y \rotatemathfig{0.04}{90}{cobordisms/neck_cutting_middle} +
 z \rotatemathfig{0.04}{90}{cobordisms/neck_cutting_right}.
\end{equation*}
We can determine the coefficients here by considering various
closures.

Adding a punctured torus at the top and a disc at the bottom gives
us $t= x t^2$, and vice versa gives us $t = z t^2$, so
$x=z=\frac{1}{t}$. Adding a `choking torus' at top and bottom gives $-t^2 = y t^4$, so $y =
-\frac{1}{t^2}$. Finally, gluing top to bottom gives $t =
\frac{1}{t} t - \frac{1}{t^2} (-t^2) + \frac{1}{t} t = 3$. We've at
this stage recovered the neck cutting relation!

\subsubsection{Consistency of the evaluation relations}%
\label{sec:evaluation-consistency}%
\begin{proof}[Proof of Lemma \ref{lem:evaluation-consistency}.]
In $\Cob{2}$, all closed foams are equivalent to scalars. This is
not as immediately apparent in $\Cob{3}$, but it's in fact true even
in $\Cob{3}^\text{ev}$; that is, even when we only impose the
evaluation relations. We describe an algorithm for evaluating closed
foams and prove that it's well-defined with respect to the
evaluation relations.

The first step, in which we do nearly all the work, is to perform
neck cutting on each sheet incident at each seam (all of which are
circles). Thus if there are $k$ seams in a closed foam, we perform
neck cutting $3k$ times, resulting in $3^{3k}$ terms. The
compensation for creating so many terms is that each term is now
relatively simple, being a disjoint union of two different types of
small closed foams.

The first type, arising from a seam in the original closed foam,
consists simply of a seam, with three of the elements appearing in Equation \eqref{eq:disc-morphisms} attached.

The second type, arising from a sheet in the original closed foam,
consists of a closed foam in which the only seams appears as part of
some `choking torus'. Notice that all of these choking toruses are
of the same type; the cyclic order around the seam is
`bulk-handle-disc', simply because this is the cyclic order appearing in the neck cutting relation. These surfaces are thus parameterized by two
numbers; the number of choking toruses, and the number of punctured
toruses. We'll write such a surface as $\Sigma_{k,l}$:
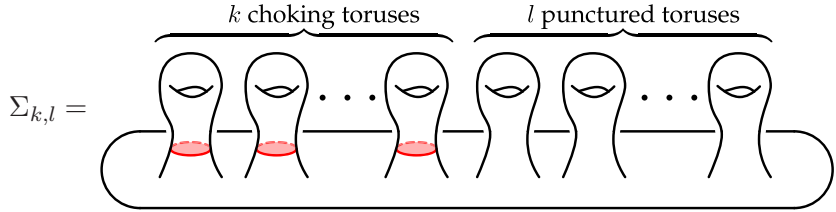
$$\Sigma_{k,l} = \mathfig{0.7}{cobordisms/sigma-kl}.$$

The second step of the algorithm is to evaluate all of these small
closed foams. In the first type, we quickly see by the seam swapping
relation that nearly all are zero. In particular, unless the three
different sheets carry different surfaces, the closed foam must be
zero. There are thus only two non-zero possibilities, depending on
the cyclic order around the seam. We can either have
`disc/handle/punctured torus' or `disc/punctured torus/handle':
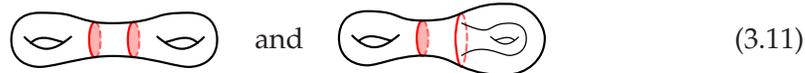
\begin{align}
\label{eq:small-closed-foams-1}
\mathfig{0.2}{cobordisms/airlock_double_torus} \quad \text{and}
\quad \mathfig{0.2}{cobordisms/airlock_double_torus_backwards}
\end{align}

We now apply the seam-swapping if we find ourselves in the second
case, then evaluate the first closed foam (via `airlock') as $-9$.

We evaluate nearly every case of the second type of closed foam,
$\Sigma_{k,l}$ by making use of Equation
\eqref{eq:torus-choking-torus}. Specifically, if $k\geq1$ and
$l\geq1$, or simply $l\geq2$, we see $\Sigma_{k,l}=0$. If $k\geq3$,
Equations \eqref{eq:choking-torus-multiplication} and
\eqref{eq:torus-choking-torus} together imply $\Sigma_{k,l}=0$. This
leaves four cases, shown in Figure
\ref{fig:irreducible-sigma-foams}, each of which we already know how
to evaluate directly.

\begin{figure}[ht]
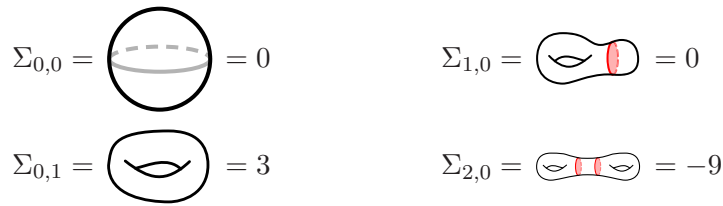

\begin{align*}
 \Sigma_{0,0} = \mathfig{0.1}{cobordisms/sphere} & = 0                 & \Sigma_{1,0} = \mathfig{0.1}{cobordisms/torus_right_blister} & = 0 \\
 \Sigma_{0,1} = \mathfig{0.1}{cobordisms/torus}  & = 3                 & \Sigma_{2,0} = \mathfig{0.1}{cobordisms/airlock_double_torus}   & = -9
\end{align*}
\caption{The irreducible examples of $\Sigma_{k,l}$, modulo neck
cutting.}\label{fig:irreducible-sigma-foams}
\end{figure}

The algorithm described so far evaluates any closed foam as a
scalar. We now check that the evaluation relations are consistent,
by showing that the evaluation algorithm produces the same result on
either side of each relation, when applied to some large closed
foam. This check requires a few cases, each of which is almost
trivial.

The first, and most trivial, cases are the closed foam relations.
It's easy to see that applying the above algorithm to any of the four closed foams in Equation \eqref{eq:closed_foams} above simply
gives the specified evaluation. This is completely trivial in 3
cases, and a short calculation for $\mathfig{0.1}{cobordisms/triple_torus_disc}$ (because there we do some
`unnecessary' neck cutting).

The seam-swapping relation is also relatively trivial. If we change
the cyclic order at a seam, the evaluation algorithm only differs in
that the two surfaces in Equation \eqref{eq:small-closed-foams-1} are
interchanged, resulting in an extra sign (actually, these two
surfaces actually occur three times each, corresponding to the three
cyclic permutations around the seam, but each pair is interchanged).

Slightly more interesting is the airlock relation. Here we simply
need to check that when we cut both seams in an airlock, modulo the
specified closed foam evaluations, we obtain exactly the other side
of the airlock relation.

Most interesting is the neck cutting relation. There are three
distinct ways we can apply the neck cutting relation; parallel to a
seam, not parallel but still separating the sheet into two pieces,
and non-separating. The first is easy; the evaluation algorithm
produces the same result, simply because neck cutting twice along
parallel circles is the same as neck cutting once (modulo evaluating
the 9 resulting closed foams). If we apply neck cutting separating a
sheet into two pieces, it's obviously the same as applying a
corresponding neck cutting to one of the second type of small closed
foams resulting from the evaluation algorithm. Thus we need to check
that the evaluation algorithm produces the same results on
$\Sigma_{k_1+k_2,l_1+l_2}$ and on
$$\frac{1}{3}\Sigma_{k_1,l_1+1}\Sigma_{k_2,l_2} - \frac{1}{9}\Sigma_{k_1+1,l_1}\Sigma_{k_2+1,l_2} + \frac{1}{3}\Sigma_{k_1,l_1}\Sigma_{k_2,l_2+1}.$$
This check involves quite a few cases; when $k_1+k_2, l_1+l_2 \geq
1$ (which splits into two subcases, $k_1,l_1 \geq 1$, and $k_1,l_2
\geq 1$), when $l_1+l_2 \geq 2$, when $k_1+k_2 \geq 3$, and the `small' cases when none
of these hold. Each case is pretty much immediate, however.

Finally, for a `non-separating' neck cutting relation we need to
check that the evaluation algorithm produces the same results on
$\Sigma_{k,l}$ ($l$ here must be at least 1) and
\begin{equation}
\label{eq:non-separating-neck-cutting-rhs}
\frac{2}{3}\Sigma_{k,l}-\frac{1}{9}\Sigma_{k+2,l-1}.
\end{equation}
If $k \geq 1$, each closed foam appearing here evaluates to $0$. If
$k=0$, everything is zero unless $l=1$, in which case the expression
in Equation \eqref{eq:non-separating-neck-cutting-rhs} is $\frac{2}{3}
3 - \frac{1}{9} (-9) = 3 = \Sigma_{0,1}$.
\end{proof}

\subsubsection{The local kernel}
\label{sec:local-kernel}%
For a given disc boundary $\bdy$ in a planar algebra $\pa$, the
`pairing tangle' has two internal discs, labeled by $\bdy$ and $\bdy^*$, with an empty external circle, and the obvious spaghetti:
$$\mathfig{0.2}{planar_tangles/pairing_tangle}$$
We'll denote the result
of inserting $x \in \pa_\bdy$ and $y \in \pa_{\bdy^*}$ simply by
$\pairing{x}{y}$.

\begin{defn}
In a spherical\footnote{A planar algebra is spherical if two planar
tangles with no points on the external disc which only differ by pulling an edge
`around the back' of the disc always act in the same way.} planar
algebra $\pa$, the `local kernel' (or maybe the `kernel of the
partition function') is the set of elements $x \in \pa_\bdy$ such that
the pairing of $x$ with any $y \in \pa_{\bdy^*}$ is zero.
\end{defn}
\begin{rem}
We need the adjective spherical here in order to give such a snappy
definition. In a possibly non-spherical planar algebra, you'd want
to say it's the set $x \in \pa_\bdy$ such that for every planar tangle
$T$, with no labels on the outer boundary and $k$ internal discs,
the first of which has label $\bdy$, and for every $k-1$ appropriate
elements of $\pa$, say $x_2, \dotsc, x_k$, the composition $T(x,
x_2, \dotsc, x_k)$ is zero.
\end{rem}

\begin{defn}
In a spherical canopolis $\mathcal{C}$, the `local kernel' is the
set of morphisms $(x:A \To B) \in \mathcal{C}_\bdy$ such that for every
$(y: C\To D) \in \mathcal{C}_{\bdy^*}$ and for any morphisms
$z:\emptyset \To \pairing{A}{C}$ and $w:\pairing{B}{D} \To
\emptyset$, the composition $w \compose \pairing{x}{y} \compose z$
is zero.
\end{defn}

It's obvious that in both cases, the local kernel is an ideal. One
can always quotient by the local kernel.

\begin{defn}
A planar algebra or canopolis is `nondegenerate' if the local kernel
is zero.
\end{defn}

\begin{lem}\label{lem:pull-to-ceiling}
Given any two webs $A$ and $B$ with common boundary $\partial^*$,
there is an isomorphism of $\Integer$-modules $$G: \Hom{\Cob{3}}{A}{B} \to
\Hom{\Cob{3}}{\emptyset}{\pairing{A^*}{B}}$$ induced by an invertible sequence of
canopolis operations.  (Here, $A^*$ denotes $A$ with its orientation reversed, so that it has boundary $\partial$.)
In particular, this isomorphism preserves membership in canopolis ideals.
\end{lem}
\begin{proof}
There is an obvious homeomorphism $h: A \cup B \cup (\partial \times [0,1]) \to \pairing{A^*}{B}$.
We define $G(F)$ to be $F$ with its boundary identification map $i$ replaced by $h \compose i$.  This
yields an isomorphism of the morphism spaces.

To see that this isomorphism is induced by canopolis operations, note that $A \cup B \cup (\partial \times [0,1])$
and $\emptyset \cup \pairing{A}{B} \cup (\emptyset \times [0,1])$ are naturally isotopic in the cylinder $D^2 \cup D^2 \cup (S^1 \times [0,1])$ (which is, of course, just a 2-sphere).
One may envision this isotopy as `pulling $A$ to the ceiling'.  Pick a nice isotopy and let $M$ denote its trace in
$\left(D^2 \cup D^2 \cup (S^1 \times [0,1]) \right) \times [0,1]$.  Because $M$ comes with an induced 2-dimensional
CW structure, it can be decomposed as a sequence $M_*$ of canopolis operations taking a foam in $\Hom{}{A}{B}$
(the inner can) to a foam in $\Hom{}{\emptyset}{\pairing{A^*}{B}}$ (the outer can).  Since $h$ is induced by the isotopy, $M_* = G$.
\end{proof}

\begin{rem}
This isomorphism does not preserve gradings of morphisms; see Lemma
$\ref{lem:su_3-shellback}$ for a statement involving gradings.
\end{rem}

\begin{cor}
In a spherical canopolis, the local kernel is generated as a canopolis ideal by the
set of morphisms $x:\emptyset \To B$ such that for every
$y: B \To \emptyset$, the composition $x \compose y$
is zero.
\end{cor}

With these definitions made, it's time to prove Lemma
\ref{lem:local-kernel-generators}.

In this section, we'll write $T = \frac{1}{2} T_\downarrow +
\frac{1}{2} T_\uparrow - T_z$ for the difference of the foams
appearing in the tube relation, and $R = R_x + R_y + R_z$ for the
sum of the foams appearing in the rocket relation. (That is, the
tube and rocket relations are $T=0$ and$R=0$.) We'll write $I$ for the canopolis ideal generated by $T$ and $R$.

\begin{proof}[Proof of Lemma \ref{lem:local-kernel-generators}.]
Let $\sum c_\alpha F_\alpha$ be an element of the local kernel of
$\Cob{3}$; that is, a linear combination of foams $F_\alpha \in \Hom{}{A}{B}$ such that every closure is zero.
By Lemma \ref{lem:pull-to-ceiling}, we may assume that $A$ is empty, and $B$ has empty boundary. We proceed
by induction on the complexity of $B$.

If $B$ is empty, then each $F_\alpha$ is equivalent to a scalar, so trivially $\sum c_\alpha
F_\alpha = 0 \in I$. If $B$ is nonempty, then an Euler characteristic
argument shows that $B$ contains a square, bigon, or circle.

Suppose $B$ contains a square.  We compose with an `identity
rocket' over the square, writing $F_\alpha = R_z \compose F_\alpha$.
Then
$$R_z \compose F_\alpha = R \compose F_\alpha - R_x \compose F_\alpha  - R_y \compose F_\alpha .$$
By definition $R \compose F_\alpha \in I$. We expand $R_x \compose F_\alpha$ as
$R_x^\text{upper} \compose R_x^\text{lower} \compose F_\alpha$, where
$$R_x^\text{lower} = \mathfig{0.10}{cobordisms/rocket_relation/rocket_x_lower} \quad \text{and} \quad
R_x^\text{upper} =
\mathfig{0.10}{cobordisms/rocket_relation/rocket_x_upper}.$$ Now $\sum R_x^\text{lower} \compose
c_\alpha F_\alpha  = R_x^\text{lower} \compose (\sum c_\alpha F_\alpha)$, and since $\sum c_\alpha F_\alpha$ is in the local kernel, so is
$R_x^\text{lower} \compose (\sum c_\alpha F_\alpha)$.  Also, $R_x^\text{lower} \compose (\sum c_\alpha F_\alpha)$ has a simpler
target than $B$, and is therefore in $I$ by our inductive
hypothesis.  Hence $R_x \compose \sum c_\alpha F_\alpha \in I$, and by
the same argument, $R_y \compose \sum c_\alpha F_\alpha \in I$.
Therefore $\sum c_\alpha F_\alpha \in I$.

The argument when $B$ contains a bigon is similar.  We express
$$F_\alpha = T_z \compose F_\alpha = \frac{1}{2} T_\downarrow \compose F_\alpha
 + \frac{1}{2} T_\uparrow \compose F_\alpha - T \compose F_\alpha .$$ By definition $T \compose F_\alpha \in I$. We write $T_\downarrow \compose F_\alpha
 = T_\downarrow^\text{upper} \compose T_\downarrow^\text{lower} \compose F_\alpha$, where
$$T_\downarrow^\text{lower} = \mathfig{0.10}{cobordisms/tube_relation/lower_bubble_kiss} \quad \text{and} \quad T_\downarrow^\text{upper} = \mathfig{0.10}{cobordisms/tube_relation/upper_bubble}.$$
$T_\downarrow^\text{lower} \compose \sum c_\alpha F_\alpha$ is in the
local kernel and has simpler target, and is therefore in $I$.  As
above, it follows simply that $T_\downarrow \compose \sum c_\alpha F_\alpha$ and $T_\uparrow \compose \sum c_\alpha F_\alpha$ are in
$I$, and therefore so is $\sum c_\alpha F_\alpha$.

Lastly, suppose $B$ contains a circle.  Then by the neck-cutting
relation, $$F_\alpha = \rotatemathfig{0.02}{90}{cobordisms/cylinder} \compose F_\alpha
= \frac{1}{3} \rotatemathfig{0.02}{90}{cobordisms/neck_cutting_left} \compose F_\alpha
      - \frac{1}{9} \rotatemathfig{0.02}{90}{cobordisms/neck_cutting_middle} \compose F_\alpha
      + \frac{1}{3} \rotatemathfig{0.02}{90}{cobordisms/neck_cutting_right} \compose F_\alpha
.$$
$\rotatemathfig{0.02}{90}{cobordisms/handle_bdy_left} \compose  \sum c_\alpha F_\alpha $ is an element
of the local kernel with simpler target, so by induction, it is in
$I$.  So $\rotatemathfig{0.02}{90}{cobordisms/neck_cutting_left} \compose  \sum c_\alpha F_\alpha \in I$. This
argument works for the other two terms in the above equation, and therefore
$\sum c_\alpha F_\alpha \in I$.
\end{proof}

\subsection{Isomorphisms}
\label{ssec:isomorphisms}%
In this section, we discover what all those local relations in
$\Cob{3}$ are really for: they imply certain isomorphisms between
objects in the category $\Mat{\Cob{3}}$. These isomorphisms should
be thought of as categorifications of relations appearing in the
$\su{3}$ spider.

Thus we set out to prove:
\begin{thm}%
\label{thm:isomorphisms}%
There are isomorphisms
\begin{align*}
 \begin{aligned}
   \mathfig{0.045}{webs/clockwise_circle} & \Iso q^{-2} \, \emptyset \directSum q^0 \, \emptyset \directSum q^2 \, \emptyset \\
   \mathfig{0.045}{webs/bubble} & \Iso q^{-1} \mathfig{0.009}{webs/tall_strand} \directSum q \mathfig{0.009}{webs/tall_strand}\\
   \mathfig{0.08}{webs/oriented_square}
     & \Iso \mathfig{0.08}{webs/two_strands_horizontal} \directSum \mathfig{0.08}{webs/two_strands_vertical}
  \end{aligned}
\end{align*}
\end{thm}
\begin{proof}
Let's define $\varphi:\mathfig{0.045}{webs/clockwise_circle} \To
q^{-2} \, \emptyset \directSum q^0 \, \emptyset \directSum
q^2 \, \emptyset$ and $\varphi^{-1} : q^{-2} \, \emptyset \directSum q^0 \, \emptyset \directSum
q^2 \, \emptyset \To
\mathfig{0.045}{webs/clockwise_circle}$ by
\begin{align*}
\varphi: & \xymatrix@C+=20mm@R+=5mm{
    & q^{-2} \, \emptyset \ar@{.}[d]|{\directSum} \\
    \mathfig{0.05}{webs/clockwise_circle}
                    \ar@{|->}[ur]_{\rotatemathfig{0.025}{90}{cobordisms/cap_bdy_left}}
                    \ar@{|->}[r]_(0.65){\frac{1}{3} \rotatemathfig{0.025}{90}{cobordisms/handle_disc_bdy_left}}
                    \ar@{|->}[dr]_{\frac{1}{3} \rotatemathfig{0.025}{90}{cobordisms/handle_bdy_left}} &
    q^0 \, \emptyset \ar@{.}[d]|{\directSum} \\
 & q^2 \, \emptyset } \\
\intertext{and}%
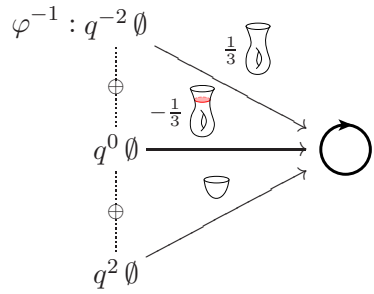
\varphi^{-1}: & \xymatrix@C+=20mm@R+=5mm{
   q^{-2} \, \emptyset \ar@{.}[d]|{\directSum} \ar[dr]^{\frac{1}{3}\rotatemathfig{0.025}{90}{cobordisms/handle_bdy_right}} & \\
    q^0 \, \emptyset \ar@{.}[d]|{\directSum} \ar[r]^(0.3){-\frac{1}{3} \rotatemathfig{0.025}{90}{cobordisms/handle_disc_bdy_right}} &
    \mathfig{0.05}{webs/clockwise_circle} \\
  q^2 \, \emptyset \ar[ur]^{\rotatemathfig{0.025}{90}{cobordisms/cap_bdy_right}} }
\end{align*}
and then perform the routine verification that these are indeed
inverses:
\begin{align*}
 \varphi^{-1} \varphi & =
    \frac{1}{3} \rotatemathfig{0.05}{90}{cobordisms/neck_cutting_right}
    - \frac{1}{9} \rotatemathfig{0.05}{90}{cobordisms/neck_cutting_middle}
    + \frac{1}{3} \rotatemathfig{0.05}{90}{cobordisms/neck_cutting_left}
    \underset{\text{neck cutting}}{=}
    \rotatemathfig{0.05}{90}{cobordisms/cylinder} = \text{id}_{\mathfig{0.025}{webs/clockwise_circle}}
\intertext{and}
 \varphi \varphi^{-1} & =
    \begin{pmatrix}
        \phantom{-\frac{1}{3}} \rotatemathfig{0.04}{90}{cobordisms/cap_bdy_right} \\
        -\frac{1}{3} \rotatemathfig{0.04}{90}{cobordisms/handle_disc_bdy_right} \\
        \phantom{-} \frac{1}{3} \rotatemathfig{0.04}{90}{cobordisms/handle_bdy_right}
    \end{pmatrix}
    \begin{pmatrix}
        \frac{1}{3} \rotatemathfig{0.04}{90}{cobordisms/handle_bdy_left} &
        \frac{1}{3} \rotatemathfig{0.04}{90}{cobordisms/handle_disc_bdy_left} &
        \rotatemathfig{0.04}{90}{cobordisms/cap_bdy_left}
    \end{pmatrix} \\
    & =
    \begin{pmatrix}
        \frac{1}{3} \mathfig{0.05}{cobordisms/torus} & \frac{1}{3} \mathfig{0.07}{cobordisms/torus_left_blister} & \mathfig{0.035}{cobordisms/sphere} \\
        -\frac{1}{9} \mathfig{0.12}{cobordisms/double_torus_disc} & -\frac{1}{9} \mathfig{0.14}{cobordisms/airlock_double_torus} & -\frac{1}{3} \mathfig{0.07}{cobordisms/torus_right_blister} \\
        \frac{1}{9} \mathfig{0.10}{cobordisms/double_torus} & \frac{1}{9} \mathfig{0.12}{cobordisms/double_torus_disc} & \frac{1}{3} \mathfig{0.05}{cobordisms/torus} \\
    \end{pmatrix} \\
    & =
    \begin{pmatrix}
       1 & 0 & 0 \\
       0 & 1 & 0 \\
       0 & 0 & 1 \\
    \end{pmatrix}
    = \text{id}_{q^{-2} \, \emptyset \directSum q^0 \, \emptyset \directSum q^2 \, \emptyset}
\end{align*}

Next we need to define the isomorphism $\mathfig{0.045}{webs/bubble}
\Iso q^{-1} \mathfig{0.009}{webs/tall_strand} \directSum q
\mathfig{0.009}{webs/tall_strand}$. It's given by
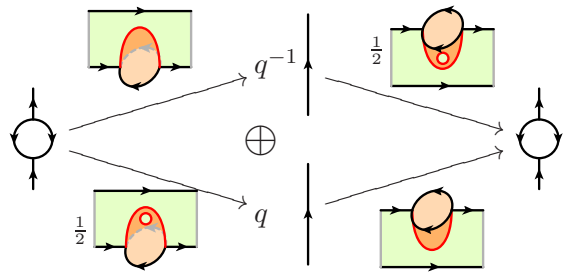
\begin{align*}
\xymatrix@C+15mm@R-15mm{%
    &
    q^{-1} \mathfig{0.009}{webs/tall_strand}
        \ar@{}[dd]_{\textrm{\scalebox{2}{$\directSum$}}}
        \ar[dr]^{\frac{1}{2}\mathfig{0.1}{cobordisms/tube_relation/upper_bubble_kiss}} &
    \\
     \mathfig{0.045}{webs/bubble}
        \ar[ur]^{\mathfig{0.1}{cobordisms/tube_relation/lower_bubble}}
        \ar[dr]_{\frac{1}{2}\mathfig{0.1}{cobordisms/tube_relation/lower_bubble_kiss}} &
    &
     \mathfig{0.045}{webs/bubble} \\
    &
    q^{\phantom{-1}} \mathfig{0.009}{webs/tall_strand}
        \ar[ur]_{\mathfig{0.1}{cobordisms/tube_relation/upper_bubble}} &
    \\
}
\end{align*}
This follows straightforwardly from the relation in Equation
\eqref{eq:tube_relation}, along with the `bagel' and `double bagel'
relations:
\begin{align}
\label{eq:bagel}%
\mathfig{0.1}{cobordisms/bagel} & = 2 &
\mathfig{0.15}{cobordisms/double_bagel} & = 0,
\end{align}
(the `bagel' here is the union of a torus and the part of the equatorial plane
outside the torus; it has two circular seams) which are easy
consequences of the `bamboo' relation appearing in Equation
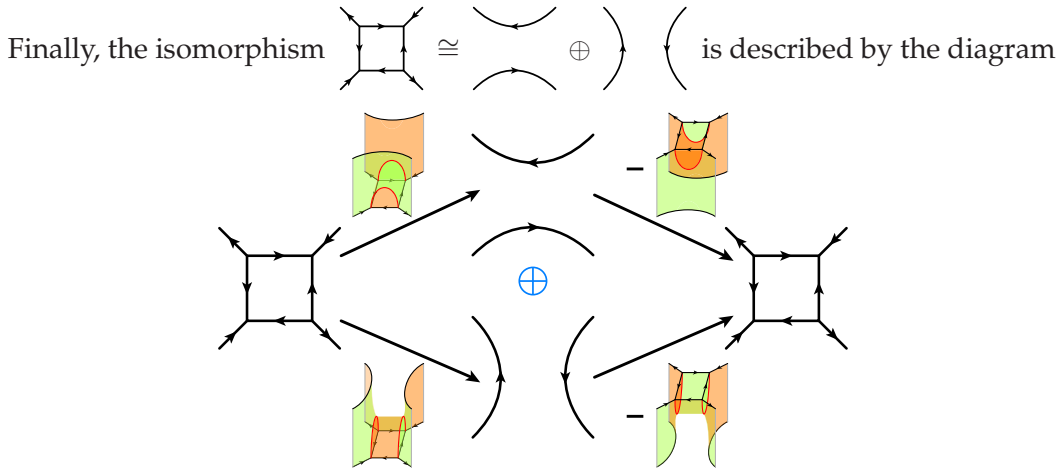
\eqref{eq:bamboo}.

Finally, the isomorphism $\mathfig{0.08}{webs/oriented_square} \Iso
\mathfig{0.08}{webs/two_strands_horizontal} \directSum
\mathfig{0.08}{webs/two_strands_vertical}$ is described by the
diagram
\begin{equation*}
\mathfig{0.6}{cobordisms/rocket_relation/isomorphism_diagram}
\end{equation*}
Verifying that these maps are mutual inverses requires the blister,
airlock and rocket relations.
\end{proof}

\section{The knot homology map}%
\label{ssec:homology_map}%
In this section we will describe the construction of the $\su{3}$
knot homology theory. This description will, of course, be
essentially equivalent to the previous constructions in
\cite{math.GT/0603307, MR2100691}, but we will emphasize certain
differences. In particular, the knot homology theory will be
explicitly local, described as a morphism of planar algebras.

The strength of this locality is that it allows us to perform
`divide and conquer' calculations. We'll explain that Bar-Natan's
\cite{math.GT/0606318} `complex simplification algorithm' can be
applied in the $\su{3}$ case. This allows us to calculate the
invariant of a knot by calculating the invariant for subtangles,
simplifying these, then gluing together the simplified complexes by
the appropriate planar operations. In \S \ref{sec:2-n-torus-knots},
we'll apply these ideas to compute the $\su{3}$ Khovanov homology of
the $(2,n)$ torus knots.

The complex simplification algorithm also allows us to give
`automatic' proofs of Reidemeister invariance; we just simplify the
complexes associated to either side of the Reidemeister move, and
observe the resulting complexes are the same.

We wish to associate to every oriented tangle a complex in
$\Mat{\Cob{3}}$.  Oriented tangles form a planar algebra generated
by the positive and negative crossings modulo relations given
by the Reidemeister moves.

In any canopolis, the complexes again form a planar algebra.
Moreover, complexes together with chain maps between them form a
canopolis. Bar-Natan proves this for $\Cob{2}$ in Theorem 2 of
\cite{MR2174270}, but his argument is completely general. There's
also a discussion of the planar algebra structure on complexes in
\cite{morrison-walker}.

It thus suffices to define the knot homology map on the positive
and negative crossing:
\begin{equation*}
\xymatrix@R-1mm{
 \mathfig{0.06}{webs/positive_crossing} \ar@{|->}[r] & & \Bigg( \bullet \ar[r]
    & q^2 \mathfig{0.06}{webs/two_strand_identity} \ar[r]^{\mathfig{0.06}{cobordisms/zip}} & q^3 \mathfig{0.06}{webs/upwards_H_diagram} \ar[r] & \bullet \Bigg) \\
 \mathfig{0.06}{webs/negative_crossing} \ar@{|->}[r] & \Bigg( \bullet \ar[r] &
    q^{-3} \mathfig{0.06}{webs/upwards_H_diagram} \ar[r]^{\mathfig{0.06}{cobordisms/unzip}} & q^{-2} \mathfig{0.06}{webs/two_strand_identity} \ar[r] & \bullet \Bigg) & \\
}
\end{equation*}

Here, the relative horizontal alignments of the complexes denote
homological height; both of the two-strand diagrams are at
homological height zero. Further, notice that, with the given
grading shifts on the objects, the differentials are grading zero
maps. Since degrees are additive under tensor products, this is true
for the differentials in the complex for any tangle.

Verifying that this map is a well-defined morphism of planar
algebras amounts to checking Reidemeister invariance, which we do in
\S \ref{ssec:isotopy-invariance}. Verifying that it's a map of
canopolises (from tangle cobordisms to chain maps) remains to be
done; we provide some evidence that this is true (on the nose, no
sign ambiguities) in \S \ref{ssec:tangle_cobordisms}.

\subsection{The simplification algorithm}
\label{ssec:simplification}%

The following lemma from \cite{math.GT/0606318} is our fundamental tool for simplifying complexes up to homotopy.

\begin{lem}[Gaussian elimination for complexes]%
\label{lem:gaussian}%
Consider the complex
\begin{equation}
\label{eq:complex}
 \xymatrix@C+55pt@R+20pt{
    A                     \ar[r]^{\psmallmatrix{\bullet \\ \alpha}}              &
    \directSumStack{B}{C} \ar[r]^{\psmallmatrix{\varphi & \lambda \\ \mu & \nu}} &
    \directSumStack{D}{E} \ar[r]^{\psmallmatrix{\bullet & \epsilon}}             &
    F
 }
\end{equation}
in any additive category, where $\varphi: B \overset{\Iso}{\To} D$
is an isomorphism, and all other morphisms are arbitrary (subject to
$d^2=0$, of course). Then there is a homotopy equivalence with a
much simpler complex, `stripping off' $\varphi$.

\begin{equation*}
 \xymatrix@C+55pt@R+20pt{
    A                     \ar[r]^{\psmallmatrix{\bullet \\ \alpha}}              \ar@{<->}[d]^{\psmallmatrix{1}} &
    \directSumStack{B}{C} \ar[r]^{\psmallmatrix{\varphi & \lambda \\ \mu & \nu}} \ar@<-0.5ex>[d]_{\psmallmatrix{0 & 1}}  &
    \directSumStack{D}{E} \ar[r]^{\psmallmatrix{\bullet & \epsilon}}             \ar@<-0.5ex>[d]_{\psmallmatrix{-\mu \varphi^{-1} & 1}} &
    F                                                                            \ar@{<->}[d]^{\psmallmatrix{1}} \\
    A \ar[r]^{\psmallmatrix{\alpha}}                          &
    C \ar[r]^{\psmallmatrix{\nu - \mu \varphi^{-1} \lambda}} \ar@<-0.5ex>[u]_{\psmallmatrix{-\varphi^{-1} \lambda \\ 1}} &
    E \ar[r]^{\psmallmatrix{\epsilon}}                       \ar@<-0.5ex>[u]_{\psmallmatrix{0 \\ 1}} &
    F
 }
\end{equation*}
\end{lem}
\begin{rem}
Note that the homotopy equivalence is also a simple homotopy
equivalence; we're just stripping off a contractible direct summand.
\end{rem}
\begin{proof}
This is simply Lemma 4.2 in \cite{math.GT/0606318} (see also Figure
2 there), this time explicitly keeping track of the chain maps.
Notice also that in a graded category, if the differentials are all
in degree $0$, so are the homotopy equivalences which we construct
here. In particular, this applies to the homotopy equivalences
associated to Reidemeister moves we construct in \S
\ref{ssec:isotopy-invariance}.
\end{proof}

We'll also state here the result of applying Gaussian elimination
twice, on two adjacent but non-composable isomorphisms. Having
these chain homotopy equivalences handy will tidy up the
calculations for the Reidemeister 2 and 3 chain maps.

\begin{lem}[Double Gaussian elimination]%
\label{lem:double_gaussian}%
When $\psi$ and $\varphi$ are isomorphisms, there's a homotopy
equivalence of complexes:
\begin{equation*}
 \xymatrix@C+55pt@R+20pt{
    A                                 \ar[r]^{\psmallmatrix{\bullet \\ \alpha}}                                    \ar@{<->}[d]^{\psmallmatrix{1}} &
    \directSumStack{B}{C}             \ar[r]^{\psmallmatrix{\psi & \beta \\ \bullet & \bullet \\ \gamma & \delta}} \ar@<-0.5ex>[d]_{\psmallmatrix{0 & 1}}  &
    \directSumStackThree{D_1}{D_2}{E} \ar[r]^{\psmallmatrix{\bullet & \varphi & \lambda \\ \bullet & \mu & \nu}}   \ar@<-0.5ex>[d]_{\psmallmatrix{-\gamma \psi^{-1} & 0 & 1}} &
    \directSumStack{F}{G}             \ar[r]^{\psmallmatrix{\bullet & \eta}}                                       \ar@<-0.5ex>[d]_{\psmallmatrix{-\mu \varphi^{-1} & 1}} &
    H                                                                                                              \ar@{<->}[d]^{\psmallmatrix{1}} \\
    A \ar[r]_{\psmallmatrix{\alpha}}                          &
    C \ar[r]_{\psmallmatrix{\delta - \gamma \psi^{-1} \beta}} \ar@<-0.5ex>[u]_<(0.25){\psmallmatrix{-\psi^{-1} \beta \\ 1}} &
    E \ar[r]_{\psmallmatrix{\nu - \mu \varphi^{-1} \lambda}}  \ar@<-0.5ex>[u]_<(0.3){\psmallmatrix{0 \\ - \varphi^{-1} \lambda \\ 1}} &
    G \ar[r]_{\psmallmatrix{\eta}}                            \ar@<-0.5ex>[u]_<(0.35){\psmallmatrix{0 \\ 1}} &
    H
 }
\end{equation*}
\end{lem}
\begin{proof}
Apply Lemma \ref{lem:gaussian} on the isomorphism $\psi$. Notice
that the isomorphism $\varphi$ survives unchanged in the resulting
complex, and apply the lemma again.
\end{proof}
\begin{rem}
Convince yourself that it doesn't matter in which order we cancel
the isomorphisms!
\end{rem}

We can now state the simplification algorithm for complexes in
$\Mat{\Cob{3}}$, analogous to Bar-Natan's algorithm
\cite{math.GT/0606318} for $\su{2}$:
\begin{itemize}
\item If an object in a complex contains a closed loop, bigon, or square, then we replace it with
the other side of the corresponding isomorphism in Theorem
\ref{thm:isomorphisms}.  (You might call this step `delooping',
`debubbling', and `desquaring'.) This increases the number of
objects in the complex, but decreases the number of possible
distinct objects, so informally we expect it to make the appearance
of isomorphisms more likely.
\item If an isomorphism appears as a matrix entry anywhere in the complex, we cancel it using Lemma \ref{lem:gaussian}.
\end{itemize}

In practice in $\Mat{\Cob{2}}$ this algorithm provides by far the
most efficient algorithm for evaluating the Khovanov homology of a
knot. This algorithm, implemented (not-so-efficiently) by Bar-Natan
and (efficiently!) by Green \cite{green-implementation} proceeds by
breaking the knot into subtangles, applying the simplification
algorithm above to the corresponding complexes, then gluing two
simplified complexes together via the appropriate planar operation,
simplifying again, and so on. Sadly, there isn't such a program for
the $\su{3}$ case.

\subsection{Isotopy invariance}
\label{ssec:isotopy-invariance}%
For each Reidemeister move, we will produce the complex associated
to the tangle on either side, and apply the simplification algorithm described above (when appropriate, also making use of Lemma
\ref{lem:double_gaussian}). There's plenty of computational work
required, but no insight. (Actually, we give a different proof of
the third Reidemeister move, trading some computation for a little
insight.) We'll produce explicit chain
maps between either side of each Reidemeister move; a gift to
whomever wants to check that the $\su{3}$ theory is functorial!

Moreover, because we use the simplification algorithm, we'll see
that the two sides of each Reidemeister move aren't just homotopic,
they're simply homotopic.\footnote{This will presumably allow an
extension of the work of Juan Ariel Ortiz-Navarro and Chris Truman
\cite{ortiz-navarro} on volume forms on Khovanov homology to the
$\su{3}$ case.}

\subsubsection{Reidemeister 1}
\label{ssec:reidemeister1}%
The complex associated to $\mathfig{0.06}{webs/positive_curl_right}$
is
$$
\xymatrix{
 q^2 \mathfig{0.06}{reidemeister_maps/R1a_complex/strand_with_circle} \ar[r]^d & q^3 \mathfig{0.06}{reidemeister_maps/R1a_complex/bigon}
}
$$
with $d$ simply a zip map. Delooping at homological height $1$, and
removing the bigon at height $2$, using the isomorphisms
\begin{align*}
 \zeta_1 & = \begin{pmatrix} \frac{1}{3} \rotatemathfig{0.025}{90}{cobordisms/handle_bdy_left}\\ \frac{1}{3} \rotatemathfig{0.025}{90}{cobordisms/handle_disc_bdy_left} \\ \rotatemathfig{0.025}{90}{cobordisms/cap_bdy_left} \end{pmatrix} &
 \zeta_2 & = \begin{pmatrix} \frac{1}{2}\mathfig{0.1}{cobordisms/tube_relation/lower_bubble_kiss} \\ \mathfig{0.1}{cobordisms/tube_relation/lower_bubble} \end{pmatrix} \\
\intertext{with inverses}
 \zeta_1^{-1} & = \begin{pmatrix} \rotatemathfig{0.025}{90}{cobordisms/cap_bdy_right} & -\frac{1}{3} \rotatemathfig{0.025}{90}{cobordisms/handle_disc_bdy_right} & \frac{1}{3}\rotatemathfig{0.025}{90}{cobordisms/handle_bdy_right} \end{pmatrix} &
 \zeta_2^{-1} & = \begin{pmatrix} \mathfig{0.1}{cobordisms/tube_relation/upper_bubble} & \frac{1}{2}\mathfig{0.1}{cobordisms/tube_relation/upper_bubble_kiss} \end{pmatrix}, \\
\end{align*}
we obtain the complex
$$
\xymatrix@C+75mm{
 \directSumStackThree{q^4 \mathfig{0.025}{reidemeister_maps/R1a_complex/strand_bending_right}}{q^2 \mathfig{0.025}{reidemeister_maps/R1a_complex/strand_bending_right}}{\phantom{q^2} \mathfig{0.025}{reidemeister_maps/R1a_complex/strand_bending_right}}
    \ar[r]^{\begin{pmatrix} \varphi = \begin{pmatrix} \mathfig{0.02}{cobordisms/sheet_algebra/sheet} & \frac{1}{6} \mathfig{0.04}{cobordisms/sheet_algebra/handle_disc} \\ 0 & -\mathfig{0.02}{cobordisms/sheet_algebra/sheet} \end{pmatrix} & \lambda = \begin{pmatrix} -\frac{1}{6} \mathfig{0.035}{cobordisms/sheet_algebra/handle} \\ \frac{1}{3} \mathfig{0.04}{cobordisms/sheet_algebra/handle_disc} \end{pmatrix}\end{pmatrix}} &
 \directSumStack{q^4 \mathfig{0.025}{reidemeister_maps/R1a_complex/strand_bending_right}}{q^2 \mathfig{0.025}{reidemeister_maps/R1a_complex/strand_bending_right}}
}.
$$

The differential here is the composition $\zeta_2 d \zeta_1^{-1}$,
and we've named some components, getting ready to apply Lemma
\ref{lem:gaussian}. Stripping off the isomorphism $\varphi$,
according to that lemma, we see that the complex is homotopy
equivalent to the desired complex: a single strand, in grading zero.
The simplifying homotopy equivalence is
\begin{align*}
 s_1 & = \begin{pmatrix}0 & 0 & \Id\end{pmatrix} \compose \zeta_1 = \mathfig{0.06}{reidemeister_maps/R1a_complex/simplifying_map} \\ s_2 & = 0 \\
\intertext{with inverse}%
 s_1^{-1} & = \zeta_1^{-1} \compose \begin{pmatrix}-\varphi^{-1} \lambda \\ \Id \end{pmatrix} = \frac{1}{3}\mathfig{0.06}{reidemeister_maps/R1a_complex/unsimplifying_map_1}-\frac{1}{9}\mathfig{0.06}{reidemeister_maps/R1a_complex/unsimplifying_map_2}+ \frac{1}{3} \mathfig{0.06}{reidemeister_maps/R1a_complex/unsimplifying_map_3} \\ s_2^{-1} & = 0.
\end{align*}

Notice here that $s_1^{-1} =
\mathfig{0.06}{reidemeister_maps/R1a_complex/unsimplifying_map_tube}$,
by the neck cutting relation. This agrees with the homotopy
equivalence proposed in \cite{math.GT/0603307}.

The calculations for the Reidemeister 1b move are much the same. We
obtain
\begin{align*}
 s_1 & =  \frac{1}{3}\mathfig{0.06}{reidemeister_maps/R1b_complex/unsimplifying_map_1}-\frac{1}{9}\mathfig{0.06}{reidemeister_maps/R1b_complex/unsimplifying_map_2}+ \frac{1}{3} \mathfig{0.06}{reidemeister_maps/R1b_complex/unsimplifying_map_3} \\
 s_2 & = 0 \\
\intertext{with inverse}%
 s_1^{-1} & = \mathfig{0.06}{reidemeister_maps/R1b_complex/simplifying_map} \\
 s_2^{-1} & = 0.
\end{align*}

\subsubsection{Reidemeister 2a}
\label{ssec:reidemeister2a}%
The complex associated to $\mathfig{0.1}{webs/R2al}$ is
$$
\xymatrix{
 q^{-1} \mathfig{0.1}{reidemeister_maps/R2a_complex/R2a_cx1} \ar[r]^{d_{-1}} &
 \directSumStack{\mathfig{0.1}{reidemeister_maps/R2a_complex/R2a_cx2a}}{\mathfig{0.1}{reidemeister_maps/R2a_complex/R2a_cx2b}} \ar[r]^{d_0} &
 q \mathfig{0.1}{reidemeister_maps/R2a_complex/R2a_cx3}
}
$$
with differentials
\begin{align*}
 d_{-1} & =
 \begin{pmatrix}
    \mathfig{0.05}{reidemeister_maps/R2a_complex/R2a_cx1_da}
    \\ \mathfig{0.05}{reidemeister_maps/R2a_complex/R2a_cx1_db}
 \end{pmatrix} &
 d_0    & =
 \begin{pmatrix}
    \mathfig{0.05}{reidemeister_maps/R2a_complex/R2a_cx2a_d}
    & - \mathfig{0.05}{reidemeister_maps/R2a_complex/R2a_cx2b_d}
 \end{pmatrix}
\end{align*}
(In this and the next section, we'll use the above shorthand for
simple foams; a red bar connecting two edges denotes a zip, and a
red bar transverse to an edge denotes an unzip.)

Applying the debubbling isomorphism $\psmallmatrix{\frac{1}{2}
\mathfig{0.05}{reidemeister_maps/R2a_complex/R2a_cx2b_debubbling_kiss}
\\ \mathfig{0.05}{reidemeister_maps/R2a_complex/R2a_cx2b_debubbling}}$ (with inverse
$\psmallmatrix{
\mathfig{0.05}{reidemeister_maps/R2a_complex/R2a_cx2b_bubbling} &
\frac{1}{2}
\mathfig{0.05}{reidemeister_maps/R2a_complex/R2a_cx2b_bubbling_kiss}}$)
to the direct summand with a bigon, we obtain the complex
$$
\xymatrix{
 q^{-1} \mathfig{0.1}{reidemeister_maps/R2a_complex/R2a_cx1} \ar[r]^{d_{-1}} &
 \directSumStackThree
    {       \mathfig{0.1}{reidemeister_maps/R2a_complex/R2a_cx2a}}
    {q      \mathfig{0.1}{reidemeister_maps/R2a_complex/R2a_cx2b_delooped}}
    {q^{-1} \mathfig{0.1}{reidemeister_maps/R2a_complex/R2a_cx2b_delooped}}
   \ar[r]^{d_0} &
 q \mathfig{0.1}{reidemeister_maps/R2a_complex/R2a_cx3}
}
$$
where
\begin{align*}
 d_{-1} & =
 \begin{pmatrix}\gamma = \mathfig{0.05}{reidemeister_maps/R2a_complex/R2a_cx1_da}
 \\ \bullet
 \\ \psi = \mathfig{0.05}{reidemeister_maps/R2a_complex/R2a_cx1}\end{pmatrix} &
 d_0    & =
 \begin{pmatrix}\lambda = \mathfig{0.05}{reidemeister_maps/R2a_complex/R2a_cx2a_d}
 & \varphi = - \mathfig{0.05}{reidemeister_maps/R2a_complex/R2a_cx2b_delooped}
 & \bullet\end{pmatrix}.
\end{align*}

Here we've named the entries of the differentials in the manner
indicated in Lemma \ref{lem:double_gaussian}. Applying that lemma
gives us chain equivalences with the desired one object complex. The
chain equivalences we're after are compositions of the chain
equivalences from Lemma \ref{lem:double_gaussian} with the debubbling
isomorphism or its inverse.

Thus the R2a `untuck' chain map is
\begin{equation*}
 \begin{pmatrix} 1 & 0 & - \gamma\psi^{-1} \end{pmatrix} \compose
 \begin{pmatrix}
    1 & 0 \\
    0 & \bullet \\
    0 & \mathfig{0.05}{reidemeister_maps/R2a_complex/R2a_cx2b_debubbling}
 \end{pmatrix}
 =
 \begin{pmatrix}
    1 & - \mathfig{0.05}{reidemeister_maps/R2a_complex/R2a_cx1_da} \compose \mathfig{0.05}{reidemeister_maps/R2a_complex/R2a_cx2b_debubbling}
 \end{pmatrix}
\end{equation*}
as claimed, and the `tuck' map is
\begin{equation*}
 \begin{pmatrix}
    1 & 0 & 0 \\
    0 & \mathfig{0.05}{reidemeister_maps/R2a_complex/R2a_cx2b_bubbling} & \bullet
 \end{pmatrix}  \compose
 \begin{pmatrix} 1 \\ - \varphi^{-1}\lambda \\ 0\end{pmatrix}
 =
 \begin{pmatrix}
    1 \\ \mathfig{0.05}{reidemeister_maps/R2a_complex/R2a_cx2b_bubbling} \compose \mathfig{0.05}{reidemeister_maps/R2a_complex/R2a_cx2a_d}
 \end{pmatrix}
\end{equation*}

\subsubsection{Reidemeister 2b}
\label{ssec:reidemeister2b}%
The complex associated to $\mathfig{0.1}{webs/R2b+}$ is
$$
\xymatrix{
 q^{-1} \mathfig{0.1}{reidemeister_maps/R2b_complex/R2b_cx1} \ar[r]^{d_{-1}} &
 \directSumStack{\mathfig{0.1}{reidemeister_maps/R2b_complex/R2b_cx2a}}{\mathfig{0.1}{reidemeister_maps/R2b_complex/R2b_cx2b}} \ar[r]^{d_0} &
 q \mathfig{0.1}{reidemeister_maps/R2b_complex/R2b_cx3}
}
$$
with differentials
\begin{align*}
 d_{-1} & =
 \begin{pmatrix}
    \mathfig{0.05}{reidemeister_maps/R2b_complex/R2b_cx1_da}
    \\ \mathfig{0.05}{reidemeister_maps/R2b_complex/R2b_cx1_db}
 \end{pmatrix} &
 d_0    & =
 \begin{pmatrix}
    - \mathfig{0.05}{reidemeister_maps/R2b_complex/R2b_cx2a_d}
    & \mathfig{0.05}{reidemeister_maps/R2b_complex/R2b_cx2b_d}
 \end{pmatrix}
\end{align*}

We now apply simplifying isomorphisms at each step (some identity sheets have been omitted in these diagrams):
\begin{align*}
\zeta_{-1} & =
\begin{pmatrix}
\frac{1}{2} \mathfig{0.1}{cobordisms/tube_relation/lower_bubble_kiss} \\
            \mathfig{0.1}{cobordisms/tube_relation/lower_bubble}
\end{pmatrix} &
\zeta_0 & =
\begin{pmatrix}
\mathfig{0.075}{cobordisms/rocket_relation/rocket_y_lower} & 0 \\
\mathfig{0.075}{cobordisms/rocket_relation/rocket_x_lower} & 0 \\
0 &             \mathfig{0.1}{reidemeister_maps/R2b_complex/R2b_cx2b_delooping} \\
0 & \frac{1}{3} \mathfig{0.1}{reidemeister_maps/R2b_complex/R2b_cx2b_delooping_handle_disc} \\
0 & \frac{1}{3} \mathfig{0.1}{reidemeister_maps/R2b_complex/R2b_cx2b_delooping_handle}
\end{pmatrix} &
\zeta_1 & =
\begin{pmatrix}
            \mathfig{0.1}{cobordisms/tube_relation/lower_bubble} \\
\frac{1}{2} \mathfig{0.1}{cobordisms/tube_relation/lower_bubble_kiss}
\end{pmatrix}
\end{align*}
with inverses (which we'll need later)
\begin{align*}
\zeta_{-1}^{-1} & =
\begin{pmatrix}
                \mathfig{0.1}{cobordisms/tube_relation/upper_bubble} &
    \frac{1}{2} \mathfig{0.1}{cobordisms/tube_relation/upper_bubble_kiss}
\end{pmatrix} &
\zeta_1^{-1} & =
\begin{pmatrix}
\frac{1}{2} \mathfig{0.1}{cobordisms/tube_relation/upper_bubble_kiss} &
            \mathfig{0.1}{cobordisms/tube_relation/upper_bubble}
\end{pmatrix}
\end{align*}
\begin{align*}
\zeta_0^{-1} & =
\begin{pmatrix}
-\mathfig{0.075}{cobordisms/rocket_relation/rocket_y_upper} & -\mathfig{0.075}{cobordisms/rocket_relation/rocket_x_upper} & 0 & 0 & 0 \\
0 & 0 &   \frac{1}{3}
\mathfig{0.1}{reidemeister_maps/R2b_complex/R2b_cx2b_looping_handle}
      & - \frac{1}{3} \mathfig{0.1}{reidemeister_maps/R2b_complex/R2b_cx2b_looping_handle_disc}
      &               \mathfig{0.1}{reidemeister_maps/R2b_complex/R2b_cx2b_looping}
\end{pmatrix}
\end{align*}

We thus obtain the complex
$$
\xymatrix{
 \directSumStack
    {q^{0}   \mathfig{0.1}{reidemeister_maps/R2b_complex/R2b_cx2b_delooped}}
    {q^{-2}  \mathfig{0.1}{reidemeister_maps/R2b_complex/R2b_cx2b_delooped}}
   \ar[r]^{d_{-1}} &
{\begin{matrix}
    q^{0}  \mathfig{0.1}{reidemeister_maps/R2b_complex/R2b_cx_simple}       \\
    \DirectSum \\
    q^{0}  \mathfig{0.1}{reidemeister_maps/R2b_complex/R2b_cx2b_delooped}   \\
    \DirectSum \\
    q^{-2} \mathfig{0.1}{reidemeister_maps/R2b_complex/R2b_cx2b_delooped}   \\
    \DirectSum \\
    q^{0}  \mathfig{0.1}{reidemeister_maps/R2b_complex/R2b_cx2b_delooped}   \\
    \DirectSum \\
    q^{+2} \mathfig{0.1}{reidemeister_maps/R2b_complex/R2b_cx2b_delooped}
 \end{matrix}}
   \ar[r]^{d_0} &
 \directSumStack
    {q^{0}   \mathfig{0.1}{reidemeister_maps/R2b_complex/R2b_cx2b_delooped}}
    {q^{+2}  \mathfig{0.1}{reidemeister_maps/R2b_complex/R2b_cx2b_delooped}}
}
$$
where
\begin{align*}
 d_{-1} & =
 \begin{pmatrix}
     \gamma = \begin{pmatrix} 0 & \mathfig{0.05}{cobordisms/saddle1} \end{pmatrix} \\
     \psi = \begin{pmatrix} - \Id & \bullet \\ 0 & \Id \end{pmatrix} \\
     \phantom{\psi = } \begin{pmatrix} \bullet & \bullet \\ \bullet & \bullet \end{pmatrix}
 \end{pmatrix} &
 d_0    & =
 \begin{pmatrix}
    \lambda = \begin{pmatrix} 0 \\ \mathfig{0.05}{cobordisms/saddle2} \end{pmatrix} & 
     \begin{pmatrix} \bullet & \bullet \\ \bullet & \bullet \end{pmatrix} &
     \varphi = \begin{pmatrix} \Id & 0 \\ \bullet & \Id \end{pmatrix}
 \end{pmatrix}.
\end{align*}

Quite a bit of cobordism arithmetic is hidden in this last step.
For example, in calculating the coefficient of the saddle appearing
$\gamma$, we used the `$\textrm{bagel} = 2$' relation. As in the R2a moves
above, we've named entries as in Lemma \ref{lem:double_gaussian},
and simply written \verb+\bullet+ for many  matrix entries,
because they won't matter in the computations to follow.

Thus the R2b `untuck' chain map is
\begin{equation*}
 \begin{pmatrix} \begin{pmatrix}1\end{pmatrix} & - \gamma\psi^{-1} = \begin{pmatrix} 0 & \mathfig{0.05}{cobordisms/saddle1} \end{pmatrix} & \begin{pmatrix} 0 & 0 \end{pmatrix} \end{pmatrix} \compose
 \begin{pmatrix}
  \mathfig{0.075}{cobordisms/rocket_relation/rocket_y_lower} & 0 \\
  \bullet & 0 \\
  0 & \mathfig{0.1}{reidemeister_maps/R2b_complex/R2b_cx2b_delooping} \\
  0 & \bullet \\
  0 & \bullet
 \end{pmatrix}
 =
 \begin{pmatrix}
    \mathfig{0.075}{cobordisms/rocket_relation/rocket_y_lower} & - \mathfig{0.075}{reidemeister_maps/R2b_complex/R2b_untuck_map_saddle}
 \end{pmatrix}
\end{equation*}
as claimed, and the `tuck' map is
\begin{equation*}
 \begin{pmatrix}
  -\mathfig{0.075}{cobordisms/rocket_relation/rocket_y_upper} & \bullet & 0 & 0 & 0 \\
  0 & 0 &   \bullet
      &  \bullet
      &               \mathfig{0.1}{reidemeister_maps/R2b_complex/R2b_cx2b_looping}
 \end{pmatrix}  \compose
 \begin{pmatrix} \phantom{- \varphi^{-1}\lambda =} \begin{pmatrix}1\end{pmatrix} \\ \phantom{- \varphi^{-1}\lambda =} \begin{pmatrix} 0 \\ 0 \end{pmatrix} \\ - \varphi^{-1}\lambda = \begin{pmatrix} 0 \\ \mathfig{0.05}{cobordisms/saddle2} \end{pmatrix} \end{pmatrix}
 =
 \begin{pmatrix}
    -\mathfig{0.075}{cobordisms/rocket_relation/rocket_y_upper} \\ - \mathfig{0.075}{reidemeister_maps/R2b_complex/R2b_tuck_map_saddle}
 \end{pmatrix}
\end{equation*}

\subsubsection{Reidemeister 3}
\label{ssec:reidemeister3} There are two almost equally appealing
approaches to the third Reidemeister move. The first is to realize
that the simplification algorithm is just as good as it is back in
the $\su{2}$ setting:
\begin{proof}[Proof modulo actually
doing all the cobordism arithmetic!]%
Apply the simplification algorithm to the complex associated to
either side of a particular variation of the third Reidemeister
move, and observe that the results are identical. Thus the two
complexes are homotopy equivalent.
\end{proof}
\begin{rem}
There's obviously some work to do here, calculating all the maps,
identifying isomorphisms, writing down the homotopy equivalences provided by Lemma \ref{lem:gaussian}, and so
on. The point is that this is all entirely algorithmic; it's an
automatic proof, with no insight required.
\end{rem}

The second method is more conceptual; it allows no real savings in
the calculations, but emphasizes that invariance
under the third Reidemeister move is a consequence of the `naturality' of the braiding in the category of complexes, described in the next two lemmas. We'll show most of the details.

\begin{lem}
Applying the simplification algorithm to the complex
\begin{equation}
\label{eq:braiding_naturality_inwards_lhs}%
\Kh{\mathfig{0.1}{webs/inwards_vertex_braiding_1}} =
 \left(
  \xymatrix@C+=20mm{
  q^4 \mathfig{0.1}{reidemeister_maps/vertex_braiding/inwards_1_cx_00}
    \ar[r]^{d_0 = \psmallmatrix{\text{zip} \\ \text{zip}}} &
  \directSumStack{q^5 \mathfig{0.1}{reidemeister_maps/vertex_braiding/inwards_1_cx_01}}{q^5 \mathfig{0.1}{reidemeister_maps/vertex_braiding/inwards_1_cx_10}}
    \ar[r]^{d_1 = \psmallmatrix{\text{zip} & -\text{zip}}} &
  q^6 \mathfig{0.1}{reidemeister_maps/vertex_braiding/inwards_1_cx_11}
  }
 \right)
\end{equation}
gives the complex
\begin{equation*}
q^8 \Kh{\mathfig{0.1}{webs/inwards_vertex_braiding_2}}\shift{+2} =
 \left(
  \xymatrix{
    q^5 \mathfig{0.1}{reidemeister_maps/vertex_braiding/inwards_2_cx_-1} \ar[r]^{\text{unzip}} &
    q^6 \mathfig{0.1}{reidemeister_maps/vertex_braiding/inwards_2_cx_0}
  }
 \right)
\end{equation*}
and the simplifying map is
\begin{align*}
 s_0 & = \begin{pmatrix} 0 \end{pmatrix} &
 s_1 & = \begin{pmatrix} 1 & - z \compose d \end{pmatrix} &
 s_2 & = \begin{pmatrix} r \end{pmatrix}.
\end{align*}
Here $d$ is the debubbling map, $z$ is a zip map, and $r$ is one of
the `half barrel' cobordisms in the `rocket isomorphism'. You can
work out exactly where all these maps are taking place simply by
considering their source and target objects.
\end{lem}
\begin{rem}
If you follow closely, you'll see we order the crossings so the
first crossing is on the right, the second crossing is on the left.
Without this, you might not like some of the signs appearing in the
proof.
\end{rem}
\begin{proof}
We begin with the complex in Equation
\eqref{eq:braiding_naturality_inwards_lhs} which, upon applying the
simplifying isomorphisms from \S \ref{ssec:isomorphisms}, becomes
$$
\xymatrix{
 q^4 \mathfig{0.1}{reidemeister_maps/vertex_braiding/inwards_1_cx_00} \ar[r]^{d_0} &
 \directSumStackThree{q^5 \mathfig{0.1}{reidemeister_maps/vertex_braiding/inwards_1_cx_01}}{q^4 \mathfig{0.1}{reidemeister_maps/vertex_braiding/inwards_1_cx_10_s}}{q^6 \mathfig{0.1}{reidemeister_maps/vertex_braiding/inwards_1_cx_10_s}} \ar[r]^{d_1} &
 \directSumStack{q^6 \mathfig{0.1}{reidemeister_maps/vertex_braiding/inwards_1_cx_11_sa}}{q^6 \mathfig{0.1}{reidemeister_maps/vertex_braiding/inwards_1_cx_11_sb}}
},
$$
with differentials
\begin{align*}
 d_0 & =
 \begin{pmatrix}
  \gamma = z \\
  \psi = \Id \\
  \bullet
 \end{pmatrix} \\
 d_1 & =
 \begin{pmatrix}
  \lambda = u & \bullet & \phi = \Id \\
  \nu = u & \bullet & \mu = 0
 \end{pmatrix},
\end{align*}
where $z$ indicates a `zip' map in the appropriate location, and $u$
an `unzip' map. Here we applied the airlock relation in calculating
$\phi$, and the blister relation in calculating $\mu$. Notice here
that $\mu = 0$, making the cancellation of the isomorphisms markedly
simple; there's no error term. We thus obtain exactly the complex
associated to $\mathfig{0.1}{webs/inwards_vertex_braiding_2}$, but
shifted up in homological height by $+2$, and in grading by $+8$.

The simplifying map itself a composition of the simplifying
isomorphisms followed by the homotopy equivalence killing off the
contractible pieces. The homotopy equivalence is $0$ at height $0$,
$\psmallmatrix{-\gamma \psi^{-1} & 0 & 1} = \psmallmatrix{-z & 0 &
1}$ at height $1$, and the identity at height $2$. Composing with
the simplifying isomorphisms gives the map in the statement of this
lemma.
\end{proof}

Analogously, we have the somewhat more awkward
\begin{lem}
The simplification algorithm provides a simple homotopy equivalence
between the complex

\begin{equation}
\label{eq:braiding_naturality_outwards_lhs}%
\Kh{\mathfig{0.1}{webs/outwards_vertex_braiding_1}} =
 \left(
  \xymatrix@C+=20mm{
  q^4 \rotatemathfig{0.1}{180}{reidemeister_maps/vertex_braiding/inwards_1_cx_00}
    \ar[r]^{d_0 = \psmallmatrix{\text{zip} \\ \text{zip}}} &
  \directSumStack{q^5 \rotatemathfig{0.1}{180}{reidemeister_maps/vertex_braiding/inwards_1_cx_10}}{q^5 \rotatemathfig{0.1}{180}{reidemeister_maps/vertex_braiding/inwards_1_cx_01}}
    \ar[r]^{d_1 = \psmallmatrix{\text{zip} & -\text{zip}}} &
  q^6 \rotatemathfig{0.1}{180}{reidemeister_maps/vertex_braiding/inwards_1_cx_11}
  }
 \right)
\end{equation}
and the complex
\begin{equation*}
 \left(
  \xymatrix{
    q^5 \rotatemathfig{0.1}{180}{reidemeister_maps/vertex_braiding/inwards_2_cx_-1} \ar[r]^{-\text{unzip}} &
    q^6 \rotatemathfig{0.1}{180}{reidemeister_maps/vertex_braiding/inwards_2_cx_0}
  }
 \right)
\end{equation*}
via the map
\begin{align*}
 s'_0 & = \begin{pmatrix} 0 \end{pmatrix} &
 s'_1 & = \begin{pmatrix} - z \compose d & 1\end{pmatrix} &
 s'_2 & = \begin{pmatrix} r \end{pmatrix}.
\end{align*}
This second complex isn't quite the complex
associated to $q^8 \mathfig{0.1}{webs/outwards_vertex_braiding_2}\shift{2}$; the
differential has been negated. Thus the map
\begin{align*}
 s''_0 & = \begin{pmatrix} 0 \end{pmatrix} &
 s''_1 & = \begin{pmatrix} z \compose d & -1 \end{pmatrix} &
 s''_2 & = \begin{pmatrix} r \end{pmatrix}.
\end{align*}
is a simple homotopy equivalence between
$q^8\Kh{\mathfig{0.1}{webs/outwards_vertex_braiding_2}}\shift{2}$ and $\Kh{\mathfig{0.1}{webs/outwards_vertex_braiding_1}}$.
\end{lem}

\begin{lem}
\label{lem:compositions-equal}%
The two compositions
\begin{align*}
\xymatrix{
 \mathfig{0.08}{reidemeister_maps/kauffman_trick/cone_source_1} \ar[r]^z &
 \mathfig{0.08}{reidemeister_maps/kauffman_trick/I_braiding_1} \ar[r]^s &
 \mathfig{0.08}{reidemeister_maps/kauffman_trick/I_braiding_2}
}%
\intertext{and}%
\xymatrix{
 \mathfig{0.08}{reidemeister_maps/kauffman_trick/cone_source_2} \ar[r]^z &
 \mathfig{0.08}{reidemeister_maps/kauffman_trick/I_braiding_3} \ar[r]^{s''} &
 \mathfig{0.08}{reidemeister_maps/kauffman_trick/I_braiding_2}
},
\end{align*}
using the maps defined in the previous two lemmas, are equal.
\end{lem}
\begin{proof}
Easy arithmetic (just in $\Integer$, not even foam arithmetic).
\end{proof}

We now need a few facts about cones.
\begin{defn}
\label{defn:cone}%
Given a chain map $f:A^\bullet \To B^\bullet$, the cone over $f$ is
$C(f)^\bullet = A^{\bullet+1} \directSum B^\bullet$, with
differential
$$d_{C(f)} = \begin{pmatrix}d_A & 0 \\ f & -d_B\end{pmatrix}$$
\end{defn}

\begin{lem}
\label{lem:postcompose_cone_morphism}
If $f:A^\bullet \To B^\bullet$ is a chain map, $r:B^\bullet \To
C^\bullet$ is a simple homotopy equivalence throwing away
contractible components (e.g. a simplification map, like those
appearing above) and $\imath:C^\bullet \To B^\bullet$ is the inverse
of $r$, then the cone $C(rf)$ is homotopic to the cone $C(f)$, via
$$\xymatrix@C+1cm{%
 C(f)^\bullet = A^{\bullet+1} \directSum B^\bullet \ar@/^/[r]^{\psmallmatrix{1 & 0 \\ 0 & r}} &
 A^{\bullet+1} \directSum C^\bullet = C(rf)^\bullet \ar@/^/[l]^{\psmallmatrix{1 & 0 \\ -hf &
 \imath}}
}$$
\end{lem}
\begin{rem}
If instead $f:B^\bullet \To A^\bullet$, then the cone $C(f\imath)$ is
homotopic to $C(f)$ via
$$\xymatrix@C+1cm{%
 C(f)^\bullet = B^{\bullet+1} \directSum A^\bullet \ar@/^/[r]^{\psmallmatrix{r & 0 \\ hf & 1}} &
 C^{\bullet+1} \directSum A^\bullet = C(f\imath)^\bullet \ar@/^/[l]^{\psmallmatrix{\imath & 0 \\ 0 & 1}}
}$$
\end{rem}

Together, the previous
four lemmas provide a proof of invariance under one variation of
the R3 move, via the categorified Kauffman trick.
\begin{align*}
 \Kh{\mathfig{0.1}{webs/R3al}} & \Iso
    \cone{\mathfig{0.08}{reidemeister_maps/kauffman_trick/cone_source_1}}{\phantom{s'' \compose}z^{\text{above}}}{\mathfig{0.08}{reidemeister_maps/kauffman_trick/I_braiding_1}} \\
        & \htpy
    \cone{\mathfig{0.08}{reidemeister_maps/kauffman_trick/cone_source_1}}{\phantom{s''}\mathllap{s} \compose z^{\text{above}}}{\mathfig{0.08}{reidemeister_maps/kauffman_trick/I_braiding_2}} \\
        & =
    \cone{\mathfig{0.08}{reidemeister_maps/kauffman_trick/cone_source_2}}{s'' \compose z^{\text{below}}}{\mathfig{0.08}{reidemeister_maps/kauffman_trick/I_braiding_2}} \\
        & \htpy
    \cone{\mathfig{0.08}{reidemeister_maps/kauffman_trick/cone_source_2}}{\phantom{s'' \compose}z^{\text{below}}}{\mathfig{0.08}{reidemeister_maps/kauffman_trick/I_braiding_3}} \\
        & \Iso \Kh{\mathfig{0.1}{webs/R3ar}}
\end{align*}

The equality on the third line is simply Lemma
\ref{lem:compositions-equal}.

The other R3 move requires similar calculations.

\subsection{Tangle cobordisms}
\label{ssec:tangle_cobordisms}%
We've almost, but not quite, provided enough detail here to check
that the $\su{3}$ cobordism theory is functorial on the nose, not
just up to sign. The calculations for the third Reidemeister move
would have to be made slightly more explicit, and then a great many
movie moves (unfortunately, there are lots of different orientations
to deal with!) need to be checked.

\begin{conj}
The $\su{3}$ cobordism theory is functorial; in particular the sign
problems seen in the $\su{2}$ case \cite{MR2174270, math.GT/0206303,
morrison-walker} don't occur.
\end{conj}
\begin{rem}
This conjecture has two sources of support. Firstly, the
representation theoretic origin of the sign problem in $\su{2}$,
namely that the standard representation is self-dual, but only
\emph{antisymmetrically} so, is simply irrelevant: the standard
representation of $\su{3}$ isn't self-dual at all. Secondly, looking
at \S \ref{ssec:reidemeister1}, we see that the coefficients of the first and last terms of
the `unsimplifying' map for the first Reidemeister move are equal.
This easily implies that the movie moves only involving the first
Reidemeister move, MM12 and MM13 (in \cite{MR2174270}'s numbering),
come out right. These moves had already failed in the $\su{2}$ case.
\end{rem}

\section{Decategorification}
\label{sec:decategorification}
\subsection{What is decategorification?}
As with quantization \cite{MR1669953}, while categorification is an
art, decategorification is a functor; it's just a fancy name for
taking the Grothendieck group\cite{wiki:Grothendieck-group}. Even
so, our situation requires slightly unusual treatment.

Usually, given an abelian category, we would form the free
$\Integer$-module on the set of objects, and add one relation $A = B
+ C$ for every short exact sequence $0 \To B \To A \To C \To 0$.

In the cobordism categories we're interested in, there are no
notions of kernels, images, or exactness.  However, our categories
still have direct sums, so we instead add relations $A = B + C$
whenever $A \Iso B \directSum C$. You can think of the result as the
`split Grothendieck group', which still makes sense in this context.

It's easy to see that we can also decategorify a canopolis; starting
with a planar algebra of categories, we obtain a planar algebra of
$\Integer$-modules.

When we decategorify a graded category, we remember the grading data and form a $\qRing$-module instead of a $\Integer$-module.

\subsection{A direct argument for $\su{2}$}
Our first result describes the decategorification of the Bar-Natan
canopolis of $\su{2}$.

\begin{defn}
The Temperley-Lieb planar algebra $\TL$ is the free planar algebra
of $\qRing$-modules with no generators, modulo the relation
$\bigcirc = q + q^{-1}$.  (Its objects are $\qRing$-linear
combinations of planar tangle diagrams modulo that relation.) The planar algebra $\TL$ is isomorphic to the representation theory of $\uqsl{2}$, or, more precisely, to the full subcategory with objects restricted to the standard representation, and tensor powers.
\end{defn}

\begin{thm} \label{thm:su_2-decategorification}
The (graded!) decategorification of the Bar-Natan canopolis
$\Mat{\Cob2}$ is the Temperley-Lieb planar algebra.
\end{thm}
\begin{proof}
The argument splits into two parts.

The first half is easy.  We must show that the relation $\bigcirc =
q + q^{-1}$ holds in the decategorification of $\Mat{\Cob2}$; that
is, $\bigcirc \Iso q \, \emptyset \directSum q^{-1} \, \emptyset$ in $\Mat{\Cob2}$. This has already been done for us by
\cite{math.GT/0606318}.

Now for the other half. We need to show that there are no more
relations in the decategorification than we one we've just seen.

Suppose we have some isomorphism $\phi : \directSum_{D} n_D D \Iso
\directSum_{D} n'_D D$, where each $D$ is a non-elliptic diagram. We
need to show that the multiplicities $n_D$ and $n'_D$ appearing on
either side agree for each diagram $D$. Fix any particular diagram
$\Delta$, let
\begin{align*}
J & = \DirectSum_{D \neq \Delta} n_D D \\
\intertext{($J$ stands for `junk'),}%
J' & = \DirectSum_{D \neq \Delta} n'_D D,
\end{align*}
and write both $\phi : n_\Delta \Delta \directSum J \To n'_\Delta
\Delta \directSum J'$ and its inverse $\phi^{-1} : n'_\Delta \Delta
\directSum J' \To n_\Delta \Delta \directSum J$ as $2 \times 2$ matrices:
\begin{align*}
\phi & = \begin{pmatrix}
 \phi_{00} : n_\Delta \Delta \To n'_\Delta \Delta & \phi_{01} : J \To n'_\Delta \Delta \\
 \phi_{10} : n_\Delta \Delta \To \mathrlap{J'}\phantom{n'_\Delta \Delta} & \phi_{11}: J \To \mathrlap{J'}\phantom{n'_\Delta \Delta}
\end{pmatrix} \\
\phi^{-1} & = \begin{pmatrix}
 \phi^{-1}_{00} : n'_\Delta \Delta \To n_\Delta \Delta & \phi^{-1}_{01} : J' \To n_\Delta \Delta \\
 \phi^{-1}_{10} : n'_\Delta \Delta \To \mathrlap{J}\phantom{n_\Delta
\Delta} & \phi^{-1}_{11}: J' \To \mathrlap{J}\phantom{n_\Delta
\Delta}
\end{pmatrix} \end{align*}

Looking at the top-left entry of the composition $\phi \phi^{-1}$,
we see that $\phi_{00} \phi^{-1}_{00} + \phi_{01} \phi^{-1}_{10}$
must be the identity on $n_\Delta \Delta$.  Notice that $\phi_{01}
\phi^{-1}_{10}$ is a linear combination of endomorphisms of
$\Delta$, each of which factors through some non-elliptic object
other than $\Delta$.  Therefore, by Corollary
\ref{cor:only-identity}, their gradings are all strictly negative,
so $\phi_{01} \phi^{-1}_{10}$ lives entirely in negative grading.
Consequently, $\phi_{00} \phi^{-1}_{00}$ is equal to the identity,
plus terms with strictly negative grading. By the same argument,
$\phi^{-1}_{00} \phi_{00}$ has the same form.  Furthermore, because
all the entries of $\phi_{00}$ and $\phi^{-1}_{00}$ are in
non-positive grading, we must have
$\left(\phi_{00}\phi^{-1}_{00}\right)_0 = \left(\phi_{00}\right)_0
\left(\phi^{-1}_{00}\right)_0$ and
$\left(\phi^{-1}_{00}\phi_{00}\right)_0 =
\left(\phi^{-1}_{00}\right)_0 \left(\phi_{00}\right)_0$. (Here the final subscript $0$ indicates the grading $0$ piece.) Therefore,
both $\left(\phi_{00}\right)_0 \left(\phi^{-1}_{00}\right)_0$ and
$\left(\phi^{-1}_{00}\right)_0 \left(\phi_{00}\right)_0$ are
identity matrices. By Corollary \ref{cor:only-identity}, the entries
of $\left(\phi_{00}\right)_0$ and $\left(\phi^{-1}_{00}\right)_0$
are simply multiples of the identity on $\Delta$.  So these two
matrices are, essentially, invertible matrices over $R$, and
therefore square \cite{wiki:Invariant-basis-number}! This gets us the desired result: $n_\Delta =
n'_\Delta$.
\end{proof}

See also \S 10 of \cite{MR2174270}, on `trace
groups', for another way to recover the Temperley-Lieb planar
algebra from this canopolis. (In fact, the construction there
doesn't really start in the same place; it uses the pure cobordism
category, whereas our decategorification only makes sense on the
category of matrices over the cobordism category, where direct sum
is defined.)

\subsection{... and why it doesn't work for $\su{3}$}
We wish to prove that the decategorification of $\Mat{\Cob3}$ is the $\su{3}$ spider: the planar algebra
of webs modulo the relations in Equation \ref{spider-relations}.

A proof along the lines of the previous section won't work for the $\su{3}$
canopolis, simply because we have no guarantee that non-identity
morphisms between non-elliptic diagrams are in negative degree. In fact, Theorem
\ref{thm:nieh-morphism} below shows that this is false.
Without this, we can't argue that (in the notation of the proof of
Theorem \ref{thm:su_2-decategorification})
$\left(\phi_{00}\phi^{-1}_{00}\right)_0 = \left(\phi_{00}\right)_0
\left(\phi^{-1}_{00}\right)_0$.

While we think it would be nice to have a proof of a $\su{3}$
decategorification statement purely in terms of the $\su{3}$
cobordism category, we'll fail at this for now, and instead describe
in \S \ref{sec:nondegeneracy} a proof that relies on some $\su{3}$
representation theory.

We'll now show that both Corollary \ref{cor:endomorphisms-grading}
and Corollary \ref{cor:only-identity} describing the morphisms in
the $\su{2}$ category fail in the $\su{3}$ category.

\begin{thm}
\label{thm:nieh-morphism}%
There are morphisms between non-elliptic objects in zero grading,
and in arbitrarily large positive gradings.
\end{thm}
\begin{proof}
See Figure \ref{fig:nieh-morphism} for the first example of a
grading zero cobordism between non-elliptic objects. We can easily
count the total grading; going from the first frame to the second,
we create $6$ circles, for a grading of $+12$, and going from the
third frame to the fourth we do $12$ `zips', for a grading of $-12$.

\begin{figure}[ht]
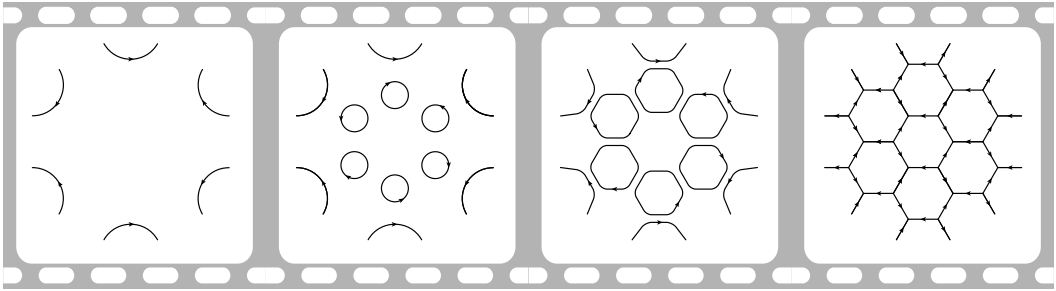

$\mathfig{1}{nieh_morphism/movie}$ \caption{The simplest example
of a grading zero cobordism between non-elliptic objects which is
not an identity cobordism.} \label{fig:nieh-morphism}
\end{figure}

Calling this cobordism $x$ and the time-reversed version $x^*$, observe that $x^* x$ is a (nonzero!)
multiple of the identity on the initial frame of Figure
\ref{fig:nieh-morphism} (and in particular, $x \neq 0$). This is an
exercise in the repeated application of the `bamboo' relation, and a
few closed foam evaluations.

We leave the construction of positive grading morphisms as an
exercise to the reader. (Hint: if you perform a sequence of zips
which produce a non-elliptic diagram with some extra circles, then
kill the resulting circles, the total grading is minus the Euler
characteristic of the graph dual to the unzipped edges.)
\end{proof}

We'll return to the consequences of this phenomenon in \S \ref{sec:karoubi}.

\subsection{Nondegeneracy}
\label{sec:nondegeneracy}
\newcommand{\emptypairing}{\pairing{\,}{\,}}
\subsubsection{Nondegeneracy for $\su{2}$}
Let $\TL_k$ denote the space of Temperley-Lieb diagrams with $k$
endpoints, modulo the usual relation $\bigcirc = q+q^{-1}$.  We
define a symmetric $\qRing$-bilinear pairing $\emptypairing_{\su{2}} :
\TL_k \times \TL_k \to \qRing$ by gluing the $k$ endpoints together, and
evaluating the resulting closed diagram.

\begin{prop} \label{prop:su_2-nondegeneracy}
The pairing $\emptypairing_{\su{2}}$ is non-degenerate on
non-elliptic diagrams.
\end{prop}
The following argument first appeared in \cite{MR1105701}.
\begin{proof}[Proof. `Diagonal dominance' \cite{wiki:Diagonally_dominant_matrix}]
Fix $k$. We'll show that the determinant of the matrix for the
pairing (with respect to the diagrammatic basis) is nonzero. This
will follow easily from the fact that the term in the determinant
corresponding to the product of the diagonal entries has strictly
higher $q$-degree than any other term.

Each entry of the matrix is of the form $(q+q^{-1})^k$, where $k$ is the
number of loops formed when two basis diagrams are glued together.
Pairing a diagram with itself produces strictly more loops than
pairing it with any other diagram, and hence the highest value of
$k$ appearing in any row appears only on the diagonal.
\end{proof}

The main result of this section is that this pairing actually tells
us the graded dimension of the space of morphisms between two
particular (unshifted) diagrams in $\Cob{2}$.

\begin{prop} \label{prop:su_2-hom pairing}
For $A$ and $B$ in $\TL_k$,
$\pairing{A}{B}_{\su{2}} = q^{\frac{k}{2}} \qdim \Hom{}{A}{B}$
\end{prop}

\begin{proof}[The easy proof specific to $\su{2}$.]

First, note that $\pairing{A}{B}_{\su{2}} = (q + q^{-1})^l$, where
$l$ is the number of boundary components of $A \cup B \cup \partial
\times [0,1]$.  By Proposition \ref{thm:su_2-classification}, the
morphism space $\Hom{}{A}{B}$ is generated by $2^l$ cobordisms consisting of $l$
connected surfaces, each of which has Euler characteristic $\pm 1$.
The degree of such a cobordism is equal to $\chi(C) - k/2$, so
$\qdim \Hom{}{A}{B} = (q + q^{-1})^l q^{-\frac{k}{2}}$, and the result
follows.
\end{proof}

However, because we have no simple classification of morphisms in $\Cob3$,
this argument does not apply to that case.  We therefore give a second
proof of Proposition \ref{prop:su_2-hom pairing}, this one using geometric
techniques that work equally well on foams.

\begin{proof}[A proof that will generalize.]

\begin{lem}[$\su{2}$ Reduction lemma]\label{lem:su_2-reduction}
Suppose $B$ contains a circle, and let $B^\bullet$ denote $B$ with
that circle removed. Then $\qdim \Hom{}{A}{B} = (q + q^{-1}) \qdim
\Hom{}{A}{B^\bullet}$, and $\pairing{A}{B}_{\su{2}} = (q + q^{-1})
\pairing{A}{B^\bullet}_{\su{2}}$.  The same result applies to
removing a circle from $A$.

\end{lem}

\begin{proof}

The first equality follows from the delooping isomorphism in
\cite{math.GT/0606318}, and the second from the definition of the
Temperley-Lieb algebra.
\end{proof}

\begin{lem}[$\su{2}$ Shellback lemma]\label{lem:su_2-shellback}
Suppose $B$ is non-elliptic and contains an arc $\alpha$ between two
adjacent boundary points.  Let $B'$ denote $B$ with $\alpha$ removed, and
let $A'$ denote $A$ with the corresponding boundary points joined by an
arc $\alpha'$.  (Note that $\partial A' = \partial B'$ has two fewer
points than $\partial A$.)  Then $\qdim \Hom{}{A}{B} = q^{-1} \qdim
\Hom{}{A'}{B'}$.
\end{lem}

\begin{proof}
Although a direct argument using canopolis operations is possible, it is
far easier to think of this operation as pulling $\alpha$ `down the
wall' of $A \cup B \cup \partial \times [0,1]$.  Because $A \cup B \cup
\partial \times [0,1]$ and $A' \cup B' \cup \partial \times [0,1]$ are
isotopic on the surface of the cylinder, there is an obvious induced
isomorphism between $\Hom{}{A}{B}$ and $\Hom{}{A'}{B'}$.  The only difference
is in the gradings, which are shifted because of the change in number of
boundary points.\end{proof}

To prove Proposition \ref{prop:su_2-hom pairing}, first observe that it
holds when $A$ and $B$ are empty diagrams.

Assume that $B$ is empty.  Since $\partial A$ is empty, $A$ is a disjoint
union of loops, and we can apply Lemma $\ref{lem:su_2-reduction}$
repeatedly to reduce to the previous case.

Assume $B$ is non-empty.  Then either $B$ contains a circle, or $B$ contains an arc
connecting adjacent boundary points.  If it contains a circle, we apply
Lemma $\ref{lem:su_2-reduction}$. Otherwise, we apply Lemma
$\ref{lem:su_2-shellback}$.  The result follows by induction on the number
of edges in $B$.
\end{proof}

We can extend this pairing to sums of diagrams:
$$\pairing{A}{B+C}_{\su{2}} = q^{\frac{k}{2}} \qdim \Hom{}{A}{B \oplus
C}.$$ (This is just observing that $\operatorname{Hom}$ respects
direct sums.)

Together, Proposition \ref{prop:su_2-nondegeneracy} and Proposition
\ref{prop:su_2-hom pairing} combine to yield a simple proof of Theorem
\ref{thm:su_2-decategorification}. Essentially, knowing that the
$\operatorname{Hom}$ pairing is nondegenerate on non-elliptic diagrams
guarantees that there are no isomorphisms amongst non-elliptic diagrams:
\begin{proof}[Alternate proof of Theorem
\ref{thm:su_2-decategorification}]
Suppose that $\oplus n_i D_i$ and $\oplus n'_i D_i$ are isomorphic
objects in $\Mat{\Cob2}$, with each $D_i$ being a non-elliptic
object. Then for any object $C$, $\qdim \Hom{}{\oplus n_i D_i}{C} =
\qdim \Hom{}{\oplus n'_i D_i}{C}$. Therefore, $$\pairing{\sum n_i D_i -
\sum n'_i D_i}{C}_{\su{2}} = 0$$ and $\sum n_i D_i = \sum n'_i D_i$
in the Temperley-Lieb algebra. There are no relations amongst
non-elliptic objects in the Temperley-Lieb planar algebra, and so
$n_i=n'_i$ for each $i$.
\end{proof}

\subsubsection{Nondegeneracy for $\su{3}$}
\label{ssec:su_3-decategorification}
We now have a new plan for a decategorification statement for
$\su{3}$; prove an analogue of Proposition \ref{prop:su_2-hom
pairing}, prove an analogue of Proposition
\ref{prop:su_2-nondegeneracy}, and then follow the alternate proof
of the $\su{2}$ decategorification statement given at the end of the
previous section.

To this end, we define a pairing $\emptypairing_{\su{3}}$ on spider
diagrams with identical boundary.  Let $\pairing{A}{B}_{\su{3}}$ be the
evaluation of the closed web resulting from reversing the orientations of
$A$, then gluing $A$ and $B$ along their boundary.  (This is $\pairing{A^*}{B}$ in the notation of \S \ref{sec:local-kernel}.)

\begin{prop} \label{prop:su_3-hom pairing}
For spider diagrams $A$ and $B$ with boundary $\partial$,
$\pairing{A}{B}_{\su{3}} = q^k \qdim \Hom{}{A}{B}$, where
$k = |\partial|$.
\end{prop}

We'll need two lemmas first. (It might be helpful to recall the isomorphisms from Theorem \ref{thm:isomorphisms} at this point.)

\begin{lem}[$\su{3}$ Reduction lemma]\label{lem:su_3-reduction}

Suppose $B$ contains a circle, and let $B^\bullet$ denote $B$ with
that circle removed.  Then $\qdim \Hom{}{A}{B} = (q^2 + 1 + q^{-2})
\qdim \Hom{}{A}{B^\bullet}$, and $\pairing{A}{B}_{\su{3}} = (q^2 + 1 +
q^{-2}) \pairing{A}{B^\bullet}_{\su{3}}$.

Similarly, assume $B$ contains a bigon, and let $B^!$ denote $B$
with that bigon deleted and replaced by an edge. Then $\qdim
\Hom{}{A}{B} = (q + q^{-1}) \qdim \Hom{}{A}{B^!}$, and
$\pairing{A}{B}_{\su{3}} = (q + q^{-1}) \pairing{A}{B^!}_{\su{3}}$.

Lastly, suppose $B$ contains a square, and let $B^\sharp$ and $B^\flat$
denote $B$ with the two possible smoothings where opposite sides of the
square are erased.  Then $\qdim \Hom{}{A}{B} = \qdim \Hom{}{A}{B^\sharp
\oplus B^\flat}$, and $\pairing{A}{B}_{\su{3}} = \pairing{A}{B^\sharp
+ B^\flat}_{\su{3}}$.

Analogous statements hold for $A$.

\end{lem}

\begin{proof}
The equalities of morphism dimensions come directly from the
isomorphisms in Theorem \ref{thm:isomorphisms}.  The equalities of
pairings are exactly Kuperberg's spider relations.
\end{proof}

\begin{lem}[$\su{3}$ Shellback lemma]\label{lem:su_3-shellback}

Suppose $B$ is non-elliptic and contains an arc $\alpha$ between two
adjacent boundary points.  Let $B'$ denote $B$ with $\alpha$ removed, and
let $A'$ denote $A$ with the corresponding boundary points joined by an
arc $\alpha'$.  Then $\qdim \Hom{}{A}{B} = q^{-2} \qdim \Hom{}{A'}{B'}$.

Suppose $B$ has a trivalent vertex $v$ with an edge $\beta$ touching
$\partial$.  Let $B^\dagger$ denote $B$ with $v$ removed and the other
edges of $v$ now terminating at $\partial$.  Let $A^\dagger$ denote $A$
with an extra vertex $v'$ added at the appropriate boundary point, and two
edges connecting it to the boundary.

Then $\qdim \Hom{}{A}{B} = q \qdim \Hom{}{A^\dagger}{B^\dagger}$.

A picture is worth far, far more than the words in the preceding
paragraph:
\begin{equation*}
 \qdim \operatorname{Hom}\left(\mathfig{0.15}{canopolis/cylinder2}\right) =
  q \qdim \operatorname{Hom}\left(\mathfig{0.15}{canopolis/cylinder1}\right)
\end{equation*}
\end{lem}

\begin{proof}
The first statement is simply Lemma \ref{lem:su_2-shellback}, modified to fit
the grading on $\su{3}$ foams.

The second looks more frightening, but it is proved by exactly the same
argument: $A \cup B \cup \partial \times [0,1]$ and $A^\dagger \cup
B^\dagger \cup \partial \times [0,1]$ are isotopic on the surface of the
cylinder, so `dragging v down the wall' changes $\Hom{}{A}{B}$ only
by a grading shift.  The power of $q$ reflects that $\partial A^\dagger$
has one more point than $\partial A$.
\end{proof}

Thus armed, we have a

\begin{proof}[Proof of Proposition \ref{prop:su_3-hom pairing}]
The proposition clearly holds when both $A$ and $B$ are empty diagrams.

Assume that $B$ is empty.  Then $A$ is a closed web, and we can
apply Lemma $\ref{lem:su_3-reduction}$ repeatedly to reduce to the
previous case.

Assume $B$ is non-empty.  If $B$ contains a circle, bigon, or square, we
apply Lemma $\ref{lem:su_3-reduction}$. Otherwise, $B$ has no closed
components, and $\partial B$ is non-empty.  In this case, either we can
find a trivalent vertex $v$ adjacent to the boundary, or $B$ is a disjoint
union of arcs, and we can find an arc $\alpha$ connecting two adjacent
boundary points.  Either one will allow us to use Lemma
$\ref{lem:su_3-shellback}$.  The result follows by induction on the number
of edges in $B$.
\end{proof}

\begin{rem}The geometrically-inclined reader may take the above nonsense with grading
shifts as evidence that a canopolis is not the most natural setting for
our seamed cobordisms.  Indeed, we claim that their native habitat is a
`spatial algebra', a higher-dimensional variant of a planar algebra.
\end{rem}

\begin{prop} \label{prop:su_3-non-degeneracy}
The pairing $\emptypairing_{\su{3}}$ is non-degenerate.
\end{prop}

It suffices to prove nondegeneracy at $q=1$, because this implies
that it holds for generic $q$.  The proof of this statement will
require an equivalent algebraic definition of
$\emptypairing_{\su{3}}$.  We can interpret any spider diagram with
boundary $\partial$ as the set of invariant tensors in $V^{\otimes
\partial}$, where $V$ is the fundamental representation of $\su{3}$. There is a standard Hermitian inner product on $V$.  If $A$ and
$B$ are spider diagrams with identical boundary, let
$\pairing{A}{B}_{R}$ denote the extension of this inner product to
tensor products of $V$ and $V^*$.  Clearly $\emptypairing_{R}$ is
nondegenerate.  It remains to show that $\emptypairing_{\su{3}}=
\emptypairing_{R}$. We will proceed, as above, by induction on $A$
and $B$.

First, if $\partial = \emptyset$, then the two pairings coincide by
\cite{MR1403861}.  For dealing with nonempty boundaries, we
prove the following lemma, which is most easily stated in pictures:

\begin{lem} \label{lem:su_3-pairing-equality}
\begin{align*}
\pairing{\mathfig{0.12}{webs_for_pairing/a_without_loop}}{\mathfig{0.12}{webs_for_pairing/b_with_loop}}_R
& =
\pairing{\mathfig{0.12}{webs_for_pairing/a_with_loop}}{\mathfig{0.12}{webs_for_pairing/b_without_loop}}_R \\
\intertext{and}
\pairing{\mathfig{0.12}{webs_for_pairing/a_without_trivalent_vertex}}{\mathfig{0.12}{webs_for_pairing/b_with_trivalent_vertex}}_R
& =
\pairing{\mathfig{0.12}{webs_for_pairing/a_with_trivalent_vertex}}{\mathfig{0.12}{webs_for_pairing/b_without_loop}}_R
\end{align*}
The corresponding statements with other orientations also hold, but we omit those calculations.
\end{lem}

Here, only the middle parts of the diagrams are meant literally; the
number of side strands is irrelevant.  In a nutshell, this says that
pieces of spider diagrams can be dragged between `floor' and
`ceiling' without changing the value of $\emptypairing_R$. Since
we know this to be the case for $\emptypairing_{\su{3}}$ by Lemma
\ref{lem:su_3-shellback}, the equality between $\emptypairing_{\su{3}}$ at $q=1$ and  $\emptypairing_R$ follows from this lemma by
induction on the size of $B$.

\begin{proof}
\newcommand{\ocup}{\mathfig{0.035}{webs_for_pairing/cup}}
\newcommand{\ocap}{\mathfig{0.035}{webs_for_pairing/cap}}
\newcommand{\trivalentCup}{\mathfig{0.05}{webs_for_pairing/trivalent_cup}}
\newcommand{\trivalentCap}{\mathfig{0.05}{webs_for_pairing/trivalent_cap}}

Translating pictures to symbols, the first statement says:
\begin{equation*}\pairing{A}{(\rm{id} \otimes \ocup \otimes \rm{id})\circ B}_R =
\pairing{(\rm{id} \otimes \ocap \otimes \rm{id})\circ A}{B}_R
\end{equation*} and the second that \begin{equation*}\pairing{A}{(\rm{id} \otimes
\trivalentCup \otimes \rm{id})\circ B}_R = \pairing{(\rm{id} \otimes
\trivalentCap \otimes \rm{id})\circ A}{B}_R\end{equation*}

Let $\{e_i\}$ be a basis for $V$ and $\{f^i\}$ the dual basis. We
write out these pictures explicitly:

\begin{align*}
\ocup & = e_1 \otimes f_1 + e_2 \otimes f_2 + e_3 \otimes f_3 \\
\ocap & = f_1 \otimes e_1 + f_2 \otimes e_2 +  f_3 \otimes e_3 \\
\trivalentCup & = \displaystyle{\sum_{\sigma \in S_3}
(-1)^{\rm{sgn}(\sigma)}} e_{\sigma(1)} \otimes e_{\sigma(2)} \otimes
e_{\sigma(3)}\\
\trivalentCap & = \displaystyle{\sum_{\sigma \in S_3}
(-1)^{\rm{sgn}(\sigma)}} f_{\sigma(1)} \otimes f_{\sigma(2)} \otimes
f_{\sigma(3)}
\end{align*}

Then the lemma follows from the definition of the inner product:
$\pairing{e_i}{e_j}_R = \delta_{ij} = \pairing{f_i}{f_j}_R$.
\end{proof}

\begin{thm} \label{thm:su_3-decategorification}
The graded decategorification of the canopolis
$\Mat{\Cob3}$ is Kuperberg's $\su{3}$ spider.
\end{thm}
\begin{rem}
See the next section, however, for a conjecture which goes further.
\end{rem}
\begin{proof}
Given Proposition \ref{prop:su_3-hom pairing}, the alternate proof of
Theorem \ref{thm:su_2-decategorification} works {\it mutatis mutandis}.
\end{proof}

\subsection{The Karoubi envelope}
\label{sec:karoubi}
We now return to the example of a degree zero non-identity morphism
from Theorem \ref{thm:nieh-morphism}. Recall we had named the
cobordism shown there in Figure \ref{fig:nieh-morphism} $x$, and
$x^*$ denoted its time reversal. We proved $x \neq 0$ by showing
$x^* x$ was a (nonzero!) multiple of the identity on the first frame.

Composing the other way round, $x x^*$ is a (multiple of a)
projection on the final frame of Figure \ref{fig:nieh-morphism}.
Normalizing correctly, let's call the projection $p$. This
projection $p$ certainly has an image in the foam category; just the
initial frame. However, $1-p$, while necessarily also being a
projection, does not have an image. (For a projection $p^2 =
p:\mathcal{O} \To \mathcal{O}$ in an arbitrary linear category, an
image is pair of morphisms $r:\mathcal{O} \To \mathcal{O'}$ and
$i:\mathcal{O'} \To \mathcal{O}$, such that $p = i \compose r$, and
$i \compose r = \Id_{\mathcal{O}'}$.) A clumsy way to see this is to
compute the pairing matrix for all non-elliptic diagrams with the
prescribed boundary; there's just a single pair of off diagonal
entries with maximal $q$ degree, corresponding via Proposition
\ref{prop:su_3-hom pairing} to the maps $r$ and $i$ for the
projection $p$, leaving no room for maps $r$ and $i$ for the
projection $1-p$.

We might suggest fixing this `problem' by passing to the Karoubi
envelope (see \cite{MR2253455} and references therein) of the foam category, which
artificially creates images for every projection. There, we can make
a conjecture relating the minimal projections appearing in the foam
category to the dual canonical basis.
\begin{conj}
The Grothendieck group of $\Kar{\Cob{3}}$ is the same as that of
$\Cob{3}$, namely the $\su{3}$ spider.

In particular, there is an ordering $\prec$ of the objects of in
$\Cob{3}$ (the ordering generated by `cap' and `unzip' will probably
do), and a bijection between non-elliptic diagrams in $\Cob{3}$ and
minimal idempotents in $\Kar{\Cob{3}}$, $D \leftrightarrow p_D$ such
that

$$\Id_D \Iso p_D \directSum \DirectSum_i q^{n_i} p_{D_i}$$

for some collection of diagrams $D_i \prec D$, and grading shifts
$n_i$. Equivalently, when we write $\Id_D$ as a sum of minimal
projections, there is one `new' projection, which we might think of
as the `leading term', plus `old' projections, each equivalent to
the new projection associated to some simpler diagram.
\end{conj}
\begin{conj}
Further, the basis for $\Kar{\Cob{3}}$ coming from the minimal
idempotents is the dual canonical basis of the $\su{3}$ spider.
\end{conj}

The immediate evidence for these conjectures is provided by the work
of Khovanov and Kuperberg in \cite{MR1680395}. There, they show that
the first non-elliptic diagram which is not a dual canonical basis
element is the final frame of the movie in Figure
\ref{fig:nieh-morphism}. Instead, in the space
$\Inv{\left(V^{\tensor2} \tensor
{V^*}^{\tensor2}\right)^{\tensor3}}$, they find that while $511$ of
the dual canonical basis vectors are given by non-elliptic diagrams,
the $512$-th is given by
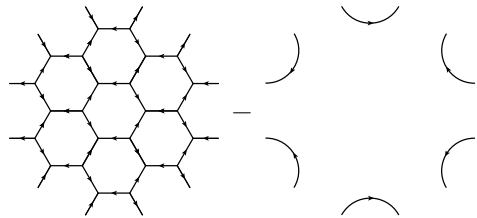
\begin{equation*}
\mathfig{0.2}{nieh_morphism/seven_hexagons} -
\mathfig{0.2}{nieh_morphism/seven_hexagons_subobject}
\end{equation*}
This is exactly the behavior described by the conjectures above. Up
until this point, every identity map on a non-elliptic diagram has
been a minimal idempotent. However, in the Karoubi envelope, we have (identifying diagrams with their identity maps)
\begin{equation*}
\mathfig{0.1}{nieh_morphism/seven_hexagons} \Iso (1-p) \directSum
\mathfig{0.1}{nieh_morphism/seven_hexagons_subobject}.
\end{equation*}

\section{Calculations}
Over $\Rational$, at least, the $\su{3}$ invariant is completely
computable for links.

\begin{lem}
For any link, there is a homotopy representation (in fact, a simple
homotopy representative) of the associated complex in the category with objects being direct sums of graded empty
diagrams and only the zero morphism.
\end{lem}
\begin{proof}
By applying isomorphisms, we can reduce the complex for a link to
one in which the objects are all direct sums of graded empty
diagrams. The morphisms are then matrices over $\Rational$; any
non-zero entry is invertible, and so there is an associated
contractible direct summand, which we can remove using Lemma \ref{lem:gaussian}.
\end{proof}

This essentially says that over $\Rational$, the homotopy type of
the invariant is characterized by its Poincar\'{e} polynomial, and
that we lose nothing by having a topological rather than
algebraic construction.

Over $\Integer[\frac{1}{2},\frac{1}{3}]$, it's more complicated; we
can still reduce all objects to the empty diagram, but there may be
`integral torsion'; the differentials may still have non-zero
entries. In the extension described in \S
\ref{sec:relaxing-relations}, in which we relax the relations
$\mathfig{0.1}{cobordisms/double_torus}  = 0$ and
$\mathfig{0.13}{cobordisms/triple_torus_disc} = 0$, there may be
further torsion associated to the polynomial ring generated by these
two foams. 

\subsection{The $(2,n)$ torus knots}
\label{sec:2-n-torus-knots}%
We now calculate the complex associated
to the two strand braid $\sigma^n$, and from that the knot homology
of its closure, the $(2,n)$ torus knot.

To begin, we introduce some notation for cobordisms,
\begin{align*}
\psi_R & = \text{unzip}_R \compose \text{zip}_R = \mathfig{0.1}{torus_knots_2n/psi_R_factor_2} \compose \mathfig{0.1}{torus_knots_2n/psi_R_factor_1} =: \mathfig{0.1}{torus_knots_2n/psi_R} \\
\psi_L & = \text{unzip}_L \compose \text{zip}_L = \mathfig{0.1}{torus_knots_2n/psi_L_factor_2} \compose \mathfig{0.1}{torus_knots_2n/psi_L_factor_1} =: \mathfig{0.1}{torus_knots_2n/psi_L}
\end{align*}
along with $\psi_\pm = \frac{1}{2}\left(\psi_R \pm \psi_L\right)$.
These cobordisms satisfy some simple relations, namely that
$\psi_R^2 = \psi_L^2 = 0$, by the double bagel relation from
Equation \eqref{eq:bagel}, and $\psi_R \psi_L = \psi_L \psi_R$. As a
consequence, $\psi_\pm \psi_\mp = 0$.

We'll further define, (harmlessly reusing names)
\begin{align*}
 \psi_R & = \mathfig{0.1}{torus_knots_2n/psi_R_bubble} &
 \psi_C & = \mathfig{0.1}{torus_knots_2n/psi_C_bubble} &
 \psi_L & = \mathfig{0.1}{torus_knots_2n/psi_L_bubble}.
\end{align*}

\newcommand{\X}{\mathfig{0.075}{torus_knots_2n/smoothing}}
\newcommand{\Y}{\mathfig{0.075}{torus_knots_2n/H}}
\newcommand{\Z}{\mathfig{0.1}{torus_knots_2n/H_with_bubble}}
\newcommand{\sX}{\mathfig{0.05}{torus_knots_2n/smoothing}}
\newcommand{\sY}{\mathfig{0.05}{torus_knots_2n/H}}
\newcommand{\sZ}{\mathfig{0.066}{torus_knots_2n/H_with_bubble}}
\newcommand{\Zreduced}{\left(\directSumStack{\Y}{\Y}\right)}

We now calculate the complex associated to a 2-twist.
\begin{thm}
Assuming $2$ is invertible, the invariant of
$\sigma^n$ is
\begin{equation}
\label{eq:2-twist-complex}
    \xymatrix{
        \X \ar[r]^{\textrm{zip}} & \Y \ar[r]^{\psi_-} & \Y \ar[r]^{\psi_+} & \cdots \ar[r]^{\psi_\mp} & \Y \ar[r]^{\psi_\pm} & \Y
    }
\end{equation}
with $\sX$ in homological height $0$, and the final $\sY$ in
homological height $n$, so the final map is $\psi_{(-1)^{n+1}}$. The
$\sX$ is in grading $2n$, the first $\sY$ in grading $2n+1$, and
each subsequent $\sY$ in grading $2$ higher than the previous, so
the last is in grading $4n-1$.
\end{thm}
\begin{proof}
The proof is by induction on $n$. For $n=1$, this complex is just
the usual invariant of a positive crossing. For $n=2$, we begin with
the complex
\newcommand{\ds}{\ar@{}[d]|<>(0.5){\scalebox{1}{$\DirectSum$}}}
\begin{equation*}
    \xymatrix@R-5mm{
        \X \ar[r]^{\textrm{zip}} \ar[dr]^{\textrm{zip}}
        & \Y \ds \ar[dr]^{\mathfig{0.05}{torus_knots_2n/psi_R_factor_1}} & \\
        & \Y \ar[r]_{-\mathfig{0.05}{torus_knots_2n/psi_L_factor_1}}
        & \Z
    }
\end{equation*}
(The sign appearing on the differential here is just the usual sign introduced by taking tensor products of complexes \cite{MR1438306}.)
Reducing the object $\sZ$ using the debubbling isomorphism, we
obtain
\begin{equation*}
    \xymatrix@R-5mm{
        \X \ar[r]^{\textrm{zip}} \ar[dr]^{\textrm{zip}}
        & \Y \ds \ar[dr]^<(0.25){\psmallmatrix{\frac{1}{2} \psi_R & 1}} & \\
        & \Y \ar[r]_<(0.25){\psmallmatrix{-\frac{\psi_L}{2} \\ -1}}
        & \Zreduced
    }
\end{equation*}
Cancelling off the matrix entry isomorphism $-1$ in the bottom row,
using Lemma \ref{lem:gaussian}, we reach the desired complex
\begin{equation*}
    \xymatrix{
        \X \ar[r]^{\textrm{zip}} & \Y \ar[r]^{\psi_-} & \Y
    }
\end{equation*}
The second differential here, $\psi_-$, is calculated as
$\frac{\psi_R}{2} - (- \frac{\psi_L}{2} \cdot (-1)^{-1} \cdot 1)$.

Now, suppose equation \eqref{eq:2-twist-complex} holds for some $n
\geq 2$. The argument is no more difficult than the $n=2$
calculation we just did, but there's more to keep track of. To
calculate $\Foam{\sigma^{n+1}}$, we simply tensor the complex in Equation \eqref{eq:2-twist-complex} with
the two step complex for a positive crossing, producing
\begin{equation*}
    \xymatrix@R-5mm{
        \X \ar[r]^{\textrm{zip}} \ar[dr]^\bullet
        & \Y \ds \ar[r]^{\psi_-} \ar[dr]^\bullet
        & \Y \ds \ar[r]^{\psi_+} \ar[dr]^\bullet
        & \cdots \ds \ar[r]^{\psi_\pm} \ar[dr]^\bullet
        & \Y \ar[dr]^{\mathfig{0.05}{torus_knots_2n/psi_R_factor_1}}   & \\
        & \Y \ar[r]_{-\mathfig{0.05}{torus_knots_2n/psi_L_factor_1}}
        & \Z \ar[r]_{- \psi_C + \psi_L}
        & \Z \ar[r]_<(0.35){- \psi_C - \psi_L}
        & \cdots \ar[r]_<(0.25){- \psi_C \mp \psi_L}
        & \Z
    }
\end{equation*}

We now reduce every diagram in the complex with the
debubbling isomorphism, obtaining
\begin{equation*}
    \xymatrix@R-5mm{
        \X \ar[r]^{\textrm{zip}} \ar[dr]^\bullet
        & \Y \ds \ar[r]^{\psi_-} \ar[dr]^\bullet
        & \Y \ds \ar[r]^{\psi_+} \ar[dr]^\bullet
        & \cdots \ds \ar[r]^{\psi_\pm} \ar[dr]^\bullet
        & \Y \ar[dr]^<(0.25){\psmallmatrix{\frac{1}{2} \psi_R & 1}}   & \\
        & \Y \ar[r]_<(0.25){\psmallmatrix{-\frac{\psi_L}{2} \\ -1}}
        & \Zreduced \ar[r]_{\Theta_+}
        & \Zreduced \ar[r]_{\Theta_-}
        & \cdots \ar[r]_{\Theta_\mp}
        & \Zreduced
    }
\end{equation*}
where $\Theta_\pm = \psmallmatrix{\pm \frac{\psi_L}{2} & 0 \\ -1 &
\pm \frac{\psi_L}{2}}$.

This complex contains many isomorphisms; we'll cancel off all the
isomorphisms appearing as matrix entries on the horizontal arrows in
the second row. This doesn't affect any of the original
differentials in the first row because there are no differentials
from the second row to the first. The
only object in the second row that survives is the first summand at
the highest homological level. The last differential is then
$\frac{1}{2} \psi_R - (\mp \frac{\psi_L}{2} \cdot (-1)^{-1} \cdot 1)
= \psi_\mp$, as claimed.

We leave it to the reader to check the gradings come out as
claimed.
\end{proof}

It's now quite easy to compute the $\su{3}$ homology
invariant for a $(2,n)$ torus knot; when we close up the braid
$\sigma^n$, all the differentials $\psi_-$ become zero, and we end
up with
\newcommand{\Xc}{\mathfig{0.075}{torus_knots_2n/smoothing_closure}}
\newcommand{\Yc}{\mathfig{0.075}{torus_knots_2n/H_closure}}
\newcommand{\psic}{\mathfig{0.05}{torus_knots_2n/psi_+_closure}}
\begin{equation*}
\begin{split}
    \xymatrix{
        q^{2n} \Xc \ar[r]^{\textrm{zip}} & q^{2n+1} \Yc \ar[r]^0 & q^{2n+3} \Yc \ar[r]^{\psic} & q^{2n+5} \Yc \ar[r]^<>(0.5)0 & \cdots
    } \\
    \xymatrix{
        \cdots \ar[r]^<>(0.5)0 & q^{4n-1} \Yc
    }
\end{split}
\end{equation*}
when $n$ is even, or
\begin{equation*}
\begin{split}
    \xymatrix{
        q^{2n} \Xc \ar[r]^{\textrm{zip}} & q^{2n+1} \Yc \ar[r]^0 & q^{2n+3} \Yc \ar[r]^{\psic} & q^{2n+5} \Yc \ar[r]^<>(0.5)0 & \cdots
    } \\
    \xymatrix{
        \cdots \ar[r]^<>(0.5){0} & q^{4n-3} \Yc \ar[r]^{\psic} & q^{4n-1} \Yc
    }
\end{split}
\end{equation*}
when $n$ is odd.

\newcommand{\Loop}{\mathfig{0.075}{torus_knots_2n/loop}}
The complex $\xymatrix{ \Xc \ar[r]^{\textrm{zip}} & q \Yc }$ is
homotopic to $\xymatrix{ q^{-2} \Loop \ar[r]^0 & \bullet }$, while
the complex $\xymatrix{ \Yc \ar[r]^{\psic} & q^2 \Yc }$ is homotopic
to $\xymatrix{ q^{-1} \Loop \ar[r]^0 & q^3 \Loop }$. Making these
replacements, we obtain the complexes
\begin{equation*}
\begin{split}
    \xymatrix{
        q^{2n-2} \Loop \ar[r]^{0} & \bullet \ar[r]^0 & q^{2n+2} \Loop \ar[r]^{0} & q^{2n+6} \Loop \ar[r]^<>(0.5)0 & \cdots
    } \\
    \xymatrix{
        \cdots \ar[r]^<>(0.5)0 & (q^{4n-2} + q^{4n}) \Loop
    }
\end{split}
\end{equation*}
when $n$ is even, or
\begin{equation*}
\begin{split}
    \xymatrix{
        q^{2n-2} \Loop \ar[r]^{0} & \bullet \ar[r]^0 & q^{2n+2} \Loop \ar[r]^{0} & q^{2n+6} \Loop \ar[r]^<>(0.5)0 & \cdots
    } \\
    \xymatrix{
        \cdots \ar[r]^<>(0.5){0} & q^{4n-4} \Loop \ar[r]^{0} & q^{4n} \Loop
    }
\end{split}
\end{equation*}
when $n$ is odd. (If you're paying careful attention to gradings, be
extra careful here; notice that the grading on the first loop
omitted by the ellipsis in the $n$ even case is actually $2n+6$
again, not $2n+10$.)

The Poincar\'{e} polynomials are thus
\begin{equation*}
  (q^{-2} + 1 + q^2)q^{2n} \left(q^{-2} + (1 + q^4 t)(q^2 t^2 + q^6 t^4 + \dotsb + q^{2n-6} t^{n-2}) + (q^{2n-2} + q^{2n})t^n\right)
\end{equation*}
when $n$ is even, and
\begin{equation*}
  (q^{-2} + 1 + q^2)q^{2n} \left(q^{-2} + (1 + q^4 t)(q^2 t^2 + q^6 t^4 + \dotsb + q^{2n-4} t^{n-1})\right)
\end{equation*}
when $n$ is odd.

The only other knot we've done calculations for is the $4_1$ knot,
whose $\su{3}$ Khovanov homology has Poincar\'{e} polynomial
$(q^{-2} + 1 + q^2)(q^{-6}t^{-2} + q^{-2}t^{-1} + 1 + q^2 t + q^6
t^2)$.

\appendix
\section{This isn't quite the same as Khovanov or Mackaay-Vaz}%
\label{sec:mackaay-vaz}

There are three significant differences between the $\su{3}$
cobordism theory defined here, and the one defined by Khovanov in
\cite{MR2100691} and deformed by Mackaay and Vaz in
\cite{math.GT/0603307}. (We assume familiarity with both of these papers throughout this section.)

The first is `locality'. Our category is described by `pictures
modulo relations', rather than by a partition function. The knot
invariant is explicitly local, defined as a map of planar algebras.

The second is that it's purely topological, in the sense that our
cobordisms don't require any dots.  As in the
$\su{2}$ case, they aren't needed, and the `sheet algebra' can be
realized by topological objects.

The third is that its deformations, in the sense of Mackaay and Vaz,
are also purely topological; instead of introducing three complex
deformation parameters, we simply remove two relations setting
certain closed foams to zero. There's a fair bit to explain here;
why, by introducing only two closed foams we see everything they see
with three deformation parameters, and the possibility of retaining
a grading in the various degenerations of the $\su{3}$ theory.

\subsection{Locality}%
Our local description of the foam category, using the canopolis
formalism, has two principal advantages over the descriptions given
in \cite{MR2100691} and \cite{math.GT/0603307}. Firstly, as
discussed previously in \S \ref{ssec:simplification}, we now have
access to Bar-Natan's simplification algorithm, which allows for
automatic proofs of Reidemeister invariance (\S \ref{ssec:isotopy-invariance}), and explicit
calculations (\S \ref{sec:2-n-torus-knots}).

Secondly, we can give a clearer analysis of the different types of
relations appearing the the theory.

Mackaay and Vaz begin by imposing certain relations on closed foams,
sufficient for evaluation; in their notation, 3D, CN, S and
$\Theta$. In our language, their Definition 2.2 says that the
category they are really interested in is the quotient by the local
kernel of the category with closed webs. (Recall the appropriate
definitions from \S \ref{sec:local-kernel}.)

Following this definition, they derive certain relations, in Lemma
2.3. We'd like to emphasize that these relations are actually of two
quite different natures. The first two, 4C (which we don't use) and
RD (our `bamboo' relation), are actually in the canopolis ideal
generated by the `evaluation' relations. On the other hand, the last
two, DR and SqR (our tube and rocket relations), cannot be derived
from the evaluation relations by canopolis operations, but only
appear in the local kernel. Moreover, while pointing out some
relations coming from the local kernel, they have no analogue of our
Lemma \ref{lem:local-kernel-generators}, providing generators of the
local kernel. Indeed, without a local setup, in which we can
describe the local kernel as a `canopolis ideal', it seems
impossible to do this.

\subsection{Relaxing our relations}
\label{sec:relaxing-relations}%
In this section, we describe a slight
generalization of our canopolis, in which we no longer impose the
relations
\begin{align*}
\mathfig{0.15}{cobordisms/double_torus} & = 0 &
\mathfig{0.2}{cobordisms/triple_torus_disc} & = 0
\end{align*}
but instead absorb these closed foams into the ground ring, calling
them $\alpha$ and $\beta$ respectively. These foams have grading
$-4$ and $-6$ respectively. This change requires modifications to
several subsequent parts of the paper.

The neck cutting relation gains an extra term\footnote{The new neck
cutting relation may be derived just as in \S
\ref{sec:explaining-relations}, being slightly more careful about
the dimensions of the various morphism spaces, taking into account
the fact that the coefficient ring is no longer all in grading
zero.}
\begin{equation}
\label{eq:neck_cutting-alpha}
 \rotatemathfig{0.05}{90}{cobordisms/cylinder} = \frac{1}{3}   \rotatemathfig{0.05}{90}{cobordisms/neck_cutting_left}
      - \frac{1}{9} \rotatemathfig{0.05}{90}{cobordisms/neck_cutting_middle}
      + \frac{1}{3} \rotatemathfig{0.05}{90}{cobordisms/neck_cutting_right}
      - \frac{1}{9} \rotatemathfig{0.05}{90}{cobordisms/neck_cutting_alpha}
\end{equation}
Consequently, there are extra terms in the sheet algebra relations,
(compare Equations \eqref{eq:choking-torus-multiplication} and
\eqref{eq:torus-choking-torus})
\begin{align*}
 \mathfig{0.075}{cobordisms/sheet_algebra/two_handle_discs} & =
  -3 \mathfig{0.075}{cobordisms/sheet_algebra/handle} - \alpha \mathfig{0.045}{cobordisms/sheet_algebra/sheet} \\
 \mathfig{0.075}{cobordisms/sheet_algebra/handle_handle_disc} & =
  \frac{2\alpha}{3}\mathfig{0.075}{cobordisms/sheet_algebra/handle_disc} + \frac{\beta}{3} \mathfig{0.045}{cobordisms/sheet_algebra/sheet} \\
 \mathfig{0.075}{cobordisms/sheet_algebra/two_handles} & =
  \frac{\alpha}{3} \mathfig{0.075}{cobordisms/sheet_algebra/handle} -
  \frac{\beta}{9} \mathfig{0.075}{cobordisms/sheet_algebra/handle_disc} + \frac{2\alpha^2}{9} \mathfig{0.045}{cobordisms/sheet_algebra/sheet}
\end{align*}
although pleasantly there are no
other changes to the local relations! These relations give the $\su{3}$ analogue of Corollary \ref{cor:su_2-sheet-algebra}.

The isomorphisms of Theorem \ref{thm:isomorphisms} mostly survive
unchanged, except the delooping isomorphism.\footnote{The authors of \cite{math.GT/0603307} don't describe a delooping isomorphism.} Now, somewhat
strangely, we have a family of isomorphisms, indexed by a parameter
$t \in S$ defined by
\begin{align*}
\varphi_t: & \xymatrix@C+=20mm@R+=5mm{
    & q^{-2} \, \emptyset \ar@{.}[d]|{\directSum} \\
    \mathfig{0.05}{webs/clockwise_circle}
                    \ar@{|->}[ur]_{\rotatemathfig{0.025}{90}{cobordisms/cap_bdy_left}}
                    \ar@{|->}[r]_(0.65){\frac{1}{3} \rotatemathfig{0.025}{90}{cobordisms/handle_disc_bdy_left}}
                    \ar@{|->}[dr]_{\frac{1}{3} \rotatemathfig{0.025}{90}{cobordisms/handle_bdy_left} + \alpha t \rotatemathfig{0.025}{90}{cobordisms/cap_bdy_left}} &
    q^0 \, \emptyset \ar@{.}[d]|{\directSum} \\
 & q^2 \, \emptyset } \\
\intertext{and}%
\varphi_t^{-1}: & \xymatrix@C+=20mm@R+=5mm{
   q^{-2} \, \emptyset \ar@{.}[d]|{\directSum} \ar[dr]^{\frac{1}{3}\rotatemathfig{0.025}{90}{cobordisms/handle_bdy_right} - \alpha (t + \frac{1}{9}) \rotatemathfig{0.025}{90}{cobordisms/cap_bdy_right}} & \\
    q^0 \, \emptyset \ar@{.}[d]|{\directSum} \ar[r]^(0.35){-\frac{1}{3} \rotatemathfig{0.025}{90}{cobordisms/handle_disc_bdy_right}} &
    \mathfig{0.05}{webs/clockwise_circle} \\
  q^2 \, \emptyset \ar[ur]^{\rotatemathfig{0.025}{90}{cobordisms/cap_bdy_right}} }
\end{align*}
It's just as easy as it was before to check that this is an
isomorphism.

Next, we turn to the isotopy invariance proofs, and check for any
use of the delooping isomorphism, or the affected relations. Both
Reidemeister 1 and Reidemeister 2b made use of the delooping
isomorphism to simplify the complexes; it turns out that the
calculation of Reidemeister 2b remains independent of which
$\varphi_t$ delooping isomorphism we use, and the chain homotopy we
produce at the end is unchanged.

The Reidemeister 1 calculation is slightly more interesting.
Using $\varphi_t$ as the delooping isomorphism, we need to modify
that calculation as follows. The isomorphisms become
\begin{align*}
 \zeta_1 & = \begin{pmatrix} \frac{1}{3} \rotatemathfig{0.025}{90}{cobordisms/handle_bdy_left} + \alpha t \rotatemathfig{0.025}{90}{cobordisms/cap_bdy_left} \\ \frac{1}{3} \rotatemathfig{0.025}{90}{cobordisms/handle_disc_bdy_left} \\ \rotatemathfig{0.025}{90}{cobordisms/cap_bdy_left} \end{pmatrix} &
 \zeta_1^{-1} & = \begin{pmatrix} \rotatemathfig{0.025}{90}{cobordisms/cap_bdy_right} & -\frac{1}{3} \rotatemathfig{0.025}{90}{cobordisms/handle_disc_bdy_right} & \frac{1}{3}\rotatemathfig{0.025}{90}{cobordisms/handle_bdy_right} - \alpha (t+\frac{1}{9}) \rotatemathfig{0.025}{90}{cobordisms/cap_bdy_right} \end{pmatrix} &
\end{align*}
and so in the differential in the simplified complex we see
$$ \lambda = \begin{pmatrix}
    -\frac{1}{6} \mathfig{0.035}{cobordisms/sheet_algebra/handle} + \alpha(t+\frac{1}{18}) \mathfig{0.02}{cobordisms/sheet_algebra/sheet} \\
    \frac{1}{3} \mathfig{0.04}{cobordisms/sheet_algebra/handle_disc}.
  \end{pmatrix}$$

Finally then, the simplifying maps acquire an extra term,
\begin{align*}
 s_1 & = \begin{pmatrix}0 & 0 & \Id\end{pmatrix} \compose \zeta_1 = \mathfig{0.06}{reidemeister_maps/R1a_complex/simplifying_map} \\ s_2 & = 0 \\
\intertext{but the inverse chain homotopy acquires an extra term}%
 s_1^{-1} & = \frac{1}{3}\mathfig{0.06}{reidemeister_maps/R1a_complex/unsimplifying_map_1}-\frac{1}{9}\mathfig{0.06}{reidemeister_maps/R1a_complex/unsimplifying_map_2}+ \frac{1}{3} \mathfig{0.06}{reidemeister_maps/R1a_complex/unsimplifying_map_3} - \frac{\alpha}{9}\mathfig{0.06}{reidemeister_maps/R1a_complex/unsimplifying_map_4} \\ s_2^{-1} & = 0.
\end{align*}

Notice, however, that it is still the case that $s_1^{-1} = \mathfig{0.06}{reidemeister_maps/R1a_complex/unsimplifying_map_tube}$.

\subsection{Dots and deformation parameters}
\newcommand{\MVF}[1]{\mathbf{Foam}_{/l}(#1)}

First, recall the definition of Mackaay and Vaz of $\MVF{a,b,c}$
(we've added the explicit notational dependence on $a$, $b$ and $c$
here).

There is an action of the ground ring on the collection of
cobordism categories $\MVF{a,b,c}$ by category equivalences:
$$\varphi_t : \mathfig{0.03}{cobordisms/sheet_algebra/sheet-1} \mapsto \mathfig{0.03}{cobordisms/sheet_algebra/sheet-1} + t \mathfig{0.03}{cobordisms/sheet_algebra/sheet},$$
taking $\MVF{a,b,c}$ to
$\MVF{a-3t,b+2at-3t^2,c+bt+at^2-t^3}$. It's easy to see that
$\varphi_{-t} \compose \varphi_t = \Id$, and that these maps
preserve the associated filtration on the categories, but not the
grading.

In the case that $a$, $b$ and $c$ are complex numbers (so the grading is already lost), it's then easy to see that $\MVF{a,b,c}$ is
isomorphic to
$\MVF{0,b+\frac{a^2}{3},c+\frac{ab}{3}+\frac{2a^3}{27}}$, and hence
we need only consider the $a=0$ case.

We now turn to showing that the dots appearing in the foams
described by Khovanov, and by Mackaay and Vaz, have `topological
representatives'. Moreover, two out of the three `deformation
parameters' in Mackaay and Vaz's paper, $b$ and $c$, also have
topological representatives.

We begin by evaluating a punctured torus in the Mackaay-Vaz theory
by neck cutting.
\begin{align}
\label{eq:MV-handle}
 \rotatemathfig{0.06}{90}{cobordisms/handle_bdy_left} & = - 3 \mathfig{0.08}{cobordisms/cap-2} + 2a\mathfig{0.08}{cobordisms/cap-1} + b\mathfig{0.08}{cobordisms/cap}
\end{align}
Next, we use the Mackaay-Vaz `bamboo' relation,\footnote{The cyclic orientation here is lower cylinder/upper cylinder/disc.}
\begin{align*}
 \rotatemathfig{0.06}{90}{cobordisms/bamboo/bamboo} = \rotatemathfig{0.06}{90}{cobordisms/bamboo/bamboo_dots1} - \rotatemathfig{0.06}{90}{cobordisms/bamboo/bamboo_dots2}
\end{align*}
to evaluate the choking torus
\begin{align}
 \label{eq:MV-handle-disc}
 \rotatemathfig{0.06}{90}{cobordisms/handle_disc_bdy_left} & = 3 \mathfig{0.08}{cobordisms/cap-1} - a \mathfig{0.08}{cobordisms/cap}
\end{align}
and thus
\begin{align*}
 \mathfig{0.08}{cobordisms/cap-1} & = 1/3 \rotatemathfig{0.06}{90}{cobordisms/handle_disc_bdy_left} + a/3 \mathfig{0.08}{cobordisms/cap}
\end{align*}

Using this, we can write any cobordism involving dots as a
$\Integer[\frac{1}{3}][a]$-linear combination of cobordisms without
dots.

What about the parameters a,b,c? Using Equation \eqref{eq:MV-handle},
we obtain
\begin{align*}
 \mathfig{0.12}{cobordisms/double_torus} & = 9 \mathfig{0.06}{cobordisms/sphere-4}  -12 a \mathfig{0.06}{cobordisms/sphere-3} + (4 a^2 - 6b) \mathfig{0.06}{cobordisms/sphere-2} \\
             & = -a^2 - 3b
\intertext{and along with Equation \eqref{eq:MV-handle-disc},}
 \mathfig{0.175}{cobordisms/triple_torus_disc} & = -a^3 - 9 ab + 27c.
\end{align*}
Rearranging these, we can express the deformation parameters $b$ and
$c$ in terms of $a$ and some closed foams.
\begin{align*}
 b & = - \frac{1}{3} \left(\mathfig{0.12}{cobordisms/double_torus} + a^2\right) \\
 c & = \frac{1}{27} \left(\mathfig{0.175}{cobordisms/triple_torus_disc} - 3a \mathfig{0.12}{cobordisms/double_torus} + a^3\right)
\end{align*}
In particular, in the special case $a=0$, we can entirely replace
the deformation parameters with closed foams.

We can now explicitly describe the correspondence between our theory
and that of Mackaay and Vaz. At the level of closed spider diagrams,
\footnote{We only say this because Mackaay and Vaz don't explicitly
describe a formalism for open diagrams.} our cobordism category is
equivalent to theirs at $a=0$, via the map
\begin{align*}
 \mathfig{0.03}{cobordisms/sheet_algebra/sheet-1} & \mapsto 1/3 \mathfig{0.06}{cobordisms/sheet_algebra/handle_disc} + a/3
     \mathfig{0.03}{cobordisms/sheet_algebra/sheet} \\
 b & \mapsto - \frac{1}{3}\mathfig{0.12}{cobordisms/double_torus} \\
 c & \mapsto \frac{1}{27} \mathfig{0.155}{cobordisms/triple_torus_disc}
\end{align*}
The inverse map is just inclusion; checking they're inverses involves a little cobordism arithmetic in each setup.

\newcommand{\urlprefix}{}
\bibliographystyle{gtart}
\bibliography{bibliography/bibliography}

This paper is available online at arXiv:\arxiv{math.GT/0612754}, and
at \url{http://scott-morrison.org/su3}.

\Addresses
\end{document}